\def\R{{\rm I}\!{\rm R}}
\def\Z{{\rm Z}\!\!{\rm Z}}
\def\N{{\rm I}\!{\rm N}}
\def\carre{\hbox{\vrule \vbox to 7pt{\hrule width 6pt \vfill \hrule}\vrule }}
\parindent = 0cm

\centerline {\bf BOUNDED WEYL PSEUDODIFFERENTIAL
OPERATORS IN FOCK SPACE}

\bigskip

\centerline {\bf (A Calder\'on-Vaillancourt theorem in an
 infinite dimensional setting)}

\bigskip

\centerline {\it Dedicated to the memory of Bernard Lascar.}

\bigskip

\centerline { L. AMOUR, L. JAGER, J. NOURRIGAT}


\bigskip

\vbox {\sevenrm
ABSTRACT 
We aim at constructing an analog of the Weyl calculus in an infinite 
dimensional setting, in which the usual configuration and phase spaces
 are 
ultimately
replaced by infinite dimensional measure spaces, the so-called abstract
 Wiener spaces. The Hilbert space on which the operators act can be seen as 
a Fock space or, equivalently, as a space of square integrable functions
on the configuration space. 
 The construction is not straightforward and needs to split 
the configuration  space 
into two factors, of which the first one is finite dimensional.
 Then one defines, for a convenient symbol $F$,
 a hybrid calculus, acting on the finite dimensional factor
 as a Weyl operator and on the other one as an anti-Wick operator, defined
thanks to an infinite dimensional Segal-Bargmann transformation.
One can establish bounds on the hybrid operators. These bounds enable us
 to prove
the convergence of  any  sequence of hybrid operators
associated with an increasing sequence of finite dimensional factors.
Their common limit is the Weyl operator $OP_h^{weyl}(F)$, the
 analog of Calder\'on-Vaillancourt Theorem being a consequence
 of the upper mentionned bounds as well.}

\bigskip
 {\sevenrm
2010  Mathematical Subject Classification  Primary 35S05; Secondary 28C20,35R15,81S30.

  }

\bigskip 
{\bf 1. Introduction.}

\bigskip

For every infinitely differentiable $F$ on $ \R^{n}\times \R^n$ with
bounded derivatives and every $h>0$, one denotes by $Op_h^{weyl
}(F)$ the bounded operator acting on  $  L^2(\R^n)$ which is
formally defined by
$$ \Big (Op_h^{weyl }(F) (\varphi ) \Big )(x) = (2 \pi h)^{-n} \int
_{\R^{2n}} e^{{i\over h} (x - y)\xi } F \left ( {x+y\over 2} , \xi
\right ) \varphi (y) dy d\xi \hskip 2cm  \varphi \in {\cal S}(\R^n).
\leqno (1.1)
$$
According to  Calder\'on-Vaillancourt Theorem [C-V], this operator
is well defined and bounded in $  L^2(\R^n)$. For further
developments on pseudo-differential operators, see for example
H\"ormander [HO] or Lerner [LER].

\bigskip

Our aim is to establish an analog of this theorem in an infinite
dimensional setting, replacing the set $\{ 1 , ... , n\}$ by a
countable set $\Gamma$, a typical example of which being a lattice
in $\R^d$.

\bigskip

The result we shall prove in this article does not exactly take the
shape of the initial theorem of [C-V] when it is restricted to a
finite dimensional space. It is rather an analog of the results
stated in  Cordes [C],
 Hwang [HW] and  Coifman Meyer [C-F].

\bigskip

Let $I_m (\{ 1 , ... , n\})$ $(m\geq 1, \ n\geq 1)$ denote the set
of  multi-indices  $(\alpha , \beta)$  such that
 $0\leq \alpha _j \leq m$ and $0 \leq \beta_j \leq m$ for all
 $j\leq n$. The results of [C], [HW] and [C-M] were not exactly concerned
with the Weyl formula (1.1), but with another one, used in those
times
 to define pseudo-differential operators. These authors prove that, if
$\partial _x^{\alpha}\partial_{\xi}^{\beta}F$ is bounded for all
 $(\alpha , \beta)$ in $I_1(\{ 1, ... , n \})$, then the  pseudo-differential
operator associated with $F$ by their formula is bounded in
 $L^2(\R^n)$.

In the case of  the Weyl calculus defined by (1.1), studying the
proof given by  A. Unterberger [U2] shows that, to ensure that
$OP_h^{weyl}(F)$ is bounded, it is sufficient to suppose that
 $\partial _x^{\alpha}\partial_{\xi}^{\beta}F$ is bounded for all
$(\alpha , \beta)$ in  $I_2(\{ 1, ... , n \})$. Our first aim is to
prove, under hypotheses similar to those of [C], [HW] and [C-M], an
upper bound on  the norm of $OP_h^{weyl}(F)$ which can easily extend
to the infinite dimensional case.

\bigskip

In the following theorem, one of the alternative statements requires
relatively
 few derivatives (at most $2$ in each  variable $x_j$ or $\xi_j$). The other
one uses derivatives up to order $4$  in each variable but yields
more precise bounds for certain examples.

\bigskip
For every non zero multi-index  $(\alpha , \beta)$, let  $S(\alpha ,
\beta)$ denote the support of  $(\alpha , \beta)$, which is the
largest subset  $S$ of $\{ 1 , ... , n\}$ such that $\alpha_j +
\beta_j \geq 1$ for all $j\in S$.

{\bf Theorem 1.1.} {\it Let $F$  be a continuous function  defined
on $\R^{2n}$. Suppose there exist  $M>0$ and a family
$(\varepsilon_j)_{(1\leq j \leq n)}$ of real nonnegative numbers
such that the following properties are satisfied

\smallskip

a) One has
$$ |F(x , \xi)| \leq M \hskip 2cm (x , \xi) \in \R^{2n} \leqno (1.2)$$
b) For every multi-index $(\alpha , \beta)$ in $I_2 ( \{1 , ... ,
n\})$, the derivative $\partial _x^{\alpha}
\partial _{\xi}^{\beta} F$ is well defined, continuous, bounded on
  $\R^{2n}$ and satisfies, if  $(\alpha , \beta)$ is not zero :
$$ |\partial _x^{\alpha} \partial _{\xi}^{\beta} F(x , \xi)| \leq M
\prod _{j\in S(\alpha , \beta)} \varepsilon _j^{\alpha_j + \beta_j}
\hskip 2cm (x , \xi) \in \R^{2n}. \leqno (1.3)$$
In this case the operator $Op_h^{weyl}(F)$ is bounded in $L^2(\R^n)$
and if $0 < h \leq 1$:
$$ \Vert Op_h^{weyl} (F) \Vert _{{\cal L} (L^2(\R^n))} \leq M \prod
_{j=1}^n (1 + 225 \pi K_2 \sqrt {h}\  \varepsilon _j) \leqno (1.4)
$$
where $K_2 = \sup _{j\leq n}\max (1 , \varepsilon_j^3)$. If
condition b) holds for every multi-index  $(\alpha , \beta)$ in
 $I_4 ( \{1 , ... , n\})$, one has :
$$ \Vert Op_h^{weyl} (F) \Vert _{{\cal L} (L^2(\R^n))} \leq M \prod
_{j=1}^n (1 + 225 \pi K_4 h \  \varepsilon _j^2 ) \leqno (1.4') $$
where $K_4= \sup _{j\leq n}\max (1 , \varepsilon_j^6)$. }

\bigskip
In order to discriminate between the roles of the $x$ and $\xi$
variables, one can alternatively assume that there exist two
sequences of nonnegative real numbers $(\rho_j)$ and $(\delta_j)$
$(1 \leq j \leq n)$ such that, for every non zero multi-index
$(\alpha , \beta)$ in $I_4 ( \{1 , ... , n\})$
$$ |\partial _x^{\alpha} \partial _{\xi}^{\beta} F(x , \xi)| \leq M
\prod _{j\in S(\alpha , \beta)} \rho _j^{\alpha_j} \delta _j^{
\beta_j} \hskip 2cm (x , \xi) \in \R^{2n} $$
In this case, one derives from the second version of Theorem 1.1
that
$$ \Vert Op_h^{weyl} (F) \Vert _{{\cal L} (L^2(\R^n))} \leq M \prod
_{j=1}^n (1 + 225 \pi K'_4 h \  \rho _j \delta_j )  $$
where $ K'_4= \sup _{j\leq n}\max (1 , (\rho_j \delta_j)^3 )$.

\bigskip

The definition of the Weyl operator and the precise  statement of
the analog of Theorem 1.1 in an infinite dimensional setting cannot
be given before Section 5. (Definition
 5.1 and Theorem  5.4),
 the necessary notions being presented
in Sections 2,3,4.
 The proof of  Theorem 1.1 is mainly contained
 in the proof of  Theorem 5.4 and will not be detailed for its own sake.

\bigskip

In the transition to an infinite dimensional situation, the set  $\{
1 , ... , n\}$ is  replaced by a countable set $\Gamma$ (for example
a lattice in
 $\R^d$ $(d\geq 1)$).
The space $L^2(\R^n)$ on which the operators act is replaced by the
 symmetrized Fock  space  ${\cal F}_s (\ell^2(\Gamma , {\bf C}))$
associated with the Hilbert space $Z =\ell^2(\Gamma , {\bf C})$.
This space  ${\cal F}_s ( \ell^2(\Gamma ,{\bf C}))$ will be denoted
by
 ${\cal H}(\Gamma)$. The definition of Fock spaces in its abstract
form is recalled in Section 3. One knows (cf. [RE-SI], [SI1], [SI2],
[J], [LEV]) that there exists an isomorphism  between
 ${\cal H}(\Gamma)$ and  $L^2 ( B(\Gamma) , \mu
_{\Gamma , h} ^K)$, where $B(\Gamma )$ is a convenient Banach space,
playing the role of the configuration space and replacing $\R^n$,
and $\mu _{\Gamma , h} ^K$ is a Gaussian measure on $B(\Gamma)$. This
isomorphism, called the Segal isomorphism, will be recalled in
Section 3. It will be denoted by  $J_{\Gamma, h}^K$.

 If the set $\Gamma$ were finite, $B(\Gamma )$ would be the space
 $ \R^{\Gamma}$ and the Gaussian measure would be defined by
$$ d\mu_{\Gamma,  h}^K (u) = (\pi h)^{-|\Gamma|/2}
 e^{-{1 \over h}|u|^2 } d\lambda_{\Gamma} (u )\ , \leqno (1.5)$$
where  $\lambda_{\Gamma}$ is the Lebesgue measure on  $\R^{\Gamma}$.
If $\Gamma$ is infinite, there is no notion of Lebesgue measure, but
the classical theory of Wiener spaces shows that the analog
 $\mu_{\Gamma, h}^K$  of the measure (1.5) can be constructed on an appropriate
Banach space  $B(\Gamma)$. This configuration space $B(\Gamma)$ is
not unique. One can choose any Banach space satisfying the
conditions required by the classical theorems related to Wiener
spaces, which can be found in
 Kuo [KU] and will be recalled in Section 2 (Theorem 2.1).
Here is an example of such a space.

\bigskip

{\bf  Definition 1.2.} {\it Let  $\Gamma$  be a countable set.
Choose a family   $b =(b_j) _{(j\in E)}$  of real positive numbers
satisfying the following property. For every $\varepsilon  >0$, the
family of positive real numbers
$$ R_j (b_j , \varepsilon ) = \int _{\varepsilon b_j}^{+\infty} e^{-{x^2 \over 2}} dx
\hskip 2cm j\in \Gamma \leqno (1.6) $$
is summable.
 Let $B_b(\Gamma)$ denote the space  of all families $(x_j) _{(j\in E)}$
such that  $\left ( {|x_j | \over b_j} \right ) _{(j\in \Gamma )}$
 converges to zero when $j$ goes to infinity. We shall
choose once for all such a family $(b_j)$ and the space
$B_b(\Gamma)$ will be denoted by $B(\Gamma)$. If $E$ is a (finite or
infinite) subset of $\Gamma$, let  $B(E)$ denote the analogous
space, corresponding to the restriction to $E$ of the same family
$(b_j)$. This space has the norm
$$ \Vert x \Vert _{B(E)} = \sup _{j\in E} {|x_j|\over b_j}\  .\leqno (1.7) $$
}

\bigskip

One shows in Section 2 (Theorem 2.3) that, for every sequence
$(b_j)_{(j\in \Gamma)}$ satisfying (1.6), the space $B_{b}(\Gamma)$
satisfies
 the hypotheses of L. Gross's Theorem 2.1 about Wiener spaces. As a
 consequence, the infinite dimensional  analog  $\mu_{\Gamma, h}^K$ of
the measure (1.5) is well defined as a measure on the Borel
$\sigma$-algebra of $B_{b}(\Gamma)$. Remark that the choice of our
configuration space does not rely on Hilbert-Schmidt operators (as
is often the case in this kind of construction), which allows us to
weaken the assumptions on the symbols. Let us give an example.

\bigskip

If $\Gamma$ is a lattice on  $\R^d$ with $d\geq 1$ and if  $| \cdot
|$ is a norm on
 $\R^d$, then, for every  $\gamma >0$, the family $b = (b_j)_{(j\in \Gamma)}$
defined by
$$ b_j = (1 + |j|)^{\gamma}$$
satisfies  the condition of Definition 1.2 (the corresponding family
$R_j (b_j, \varepsilon)$ is summable for every $\varepsilon >0$).

\bigskip

In the infinite dimensional version of Theorem 1.1, the functions
$F$ (the symbols of the operators) are bounded and continuous on the
phase space corresponding to the set $B(\Gamma)$ of Definition 1.2,
that is $B(\Gamma) \times B(\Gamma)$. We can now list the hypothesis
on the partial derivatives of the symbols. A multi-index is an
application $\alpha$ from  $\Gamma$ in $\N$, such that  $\alpha _j =
0$ except for a finite number of indices $j\in \Gamma$. We denote by
$I_m (\Gamma)$ $(m\geq 1)$ the set of multi-indices  $(\alpha ,
\beta)$ such that  $0 \leq \alpha _j \leq m$ and $0 \leq  \beta _j
\leq m $ for all  $j$ in $\Gamma$. For every non zero multi-index
$(\alpha , \beta)$, $S(\alpha , \beta)$ denotes the largest set $S$
such that $\alpha _j + \beta_j \geq 1$ for all $j$ in $S$ and is
called
 the {\it support} of  $(\alpha , \beta)$. It is therefore finite.

\bigskip

{\bf  Definition 1.3.} {\it Let $\Gamma$ be an infinite, countable
set. Let $B(\Gamma)$ be the space of Definition 1.2. Let
$\varepsilon = (\varepsilon_j )_{(j\in \Gamma)}$ be a family of
nonnegative real numbers, indexed by the elements of $\Gamma$  and
let $M$ be positive. A function $F$, bounded and continuous on
$B(\Gamma)\times B(\Gamma)$, is said to satisfy the hypothese
$H_m(M, \varepsilon)$ (with $m\geq 1)$ if
\smallskip
a)For all $(x , \xi)$ in $B(\Gamma )\times B(\Gamma)$:
$$ |F(x , \xi)| \leq M\ ; \leqno (1.7) $$
b) For every  multi-index  $(\alpha , \beta)$ in $I_m(\Gamma)$, the
partial derivative
 $\partial_x^{\alpha}\partial_{\xi}^{\beta}F$
is well defined, continuous and bounded on $B(\Gamma )\times
B(\Gamma)$ and satisfies, for all $(x , \xi)$ in $B(\Gamma )\times
B(\Gamma)$
$$ |\partial_x^{\alpha}\partial_{\xi}^{\beta}F (x , \xi)| \leq M
\prod _{j\in S(\alpha , \beta)} \varepsilon_j ^{\alpha_j + \beta_j}\
. \leqno (1.8)$$ }

\bigskip

In addition to the Weyl calculus we shall need the anti-Wick
calculus and consequently the Segal Bargmann transform in an
infinite dimensional setting. One first basically defines this
transform (Definition 4.1) as a partial isometry from any Fock space
$ {\cal F}_s ( Z)$ associated with a separable Hilbert space $Z$,
into the Fock space
 ${\cal F}_s (Z \times Z)$. In the case when
 $Z_{\bf C} =  \ell^2(\Gamma , {\bf C})$, the Segal Bargmann transform,
denoted by  $W_{\Gamma}$, is a partial isometry from the
``configuration'' Fock space  ${\cal F}_s ( \ell^2(\Gamma , {\bf
C}))= {\cal H}(\Gamma)$ (on which the operators are defined) into
the ``phase'' Fock space  ${\cal F}_s ( \ell^2(\Gamma , {\bf C})
\times \ell^2(\Gamma , {\bf C}))$, which will be denoted by ${\cal
H}_{\Phi}(\Gamma)$. This space is isomorphic to $L^2 (B(\Gamma)
\times B(\Gamma), \mu _{\Gamma , h}^{\Phi})$, where $\mu _{\Gamma ,
h}^{\Phi}$ is a Gaussian measure on
 $B(E) \times B(E)$. This Segal isomorphism is denoted by
 $J_{\Gamma , h}^{\Phi}$.
Composing the isometry $W_{\Gamma}$ with both Segal isomorphisms
yields a partial isometry from  $L^2(B(\Gamma), \mu_{\Gamma ,h}^K)$
into
 $L^2(B(\Gamma) \times B(\Gamma), \mu_{\Gamma ,h}^{\Phi})$. We shall see
in Section 4 the relation with the Bargmann transform  as it is
defined in  Kree R\c aczka [K-R]. Moreover we shall give some
equivalent characterizations of the subspace $SB(E , h)$
 of  $L^2(B(\Gamma) \times B(\Gamma), \mu_{\Gamma ,h}^{\Phi})$, which is
the range of $L^2(B(\Gamma), \mu_{\Gamma ,h}^K)$  by this
application. One of these properties is due to Driver-Hall [D-H].

\bigskip

The composition  $J_{\Gamma  h}^{\Phi} \circ W_{\Gamma}$ is a
partial isometry from  ${\cal H}(\Gamma)$ into $L^2 (B(\Gamma)
\times B(\Gamma), \mu _{\Gamma , h}^{\Phi})$. Let $F$ be bounded on
$B(\Gamma) \times B(\Gamma)$ and measurable with
 respect to the Borel $\sigma$-algebra. For all $h>0$, one can associate
 with $F$ an Anti Wick operator   $Op_h^{AW}(F)$, which is bounded in
 ${\cal H}(\Gamma)$ and satisfies, for all
$f$ and $g$ in ${\cal H}(\Gamma)$:
$$ < Op_h^{AW} (F) f , g> _{{\cal H}(\Gamma)} = \int _{B(\Gamma) \times
B(\Gamma) } F(X) \ \Big ( J_{\Gamma , h}^{\Phi} W_{\Gamma}f \Big )
(X) \overline {\Big ( J_{\Gamma , h}^{\Phi} W_{\Gamma}g \Big ) (X)}
d\mu_{\Gamma , h}^{\Phi} (X) \ . \leqno (1.)
$$

\bigskip

In this article we shall associate an operator depending on $h>0$
with every function $F$ continuous and bounded on $B(\Gamma) \times
B(\Gamma)$ and satisfying

\smallskip

- either the hypothesis  $H_2(M, \varepsilon)$, where $M>0$ and
$\varepsilon = (\varepsilon_j )_{(j\in \Gamma)}$ is a summable
family of real nonnegative numbers;

\smallskip

- or  the hypothesis  $H_4(M, \varepsilon)$, where the family
$(\varepsilon_j ^2)_{(j\in \Gamma)}$ is summable.

\smallskip

In either case the associated operator $Op_h^{weyl}(F)$ will be
bounded in the Fock space
 ${\cal H}(\Gamma) = {\cal F}_s
(\ell^2(\Gamma , {\bf C}))$.

\bigskip

The precise definition of this operator cannot be given yet;
 let us just say that, for every
finite subset $E$ of $\Gamma$, one defines a hybrid operator
 $Op_h^{hyb , E}(F)$, acting as a Weyl operator with respect to the
variables $x_j$ $(j\in E)$ and as an Anti Wick operator with respect
to the  variables $x_k$ $(k \notin E)$. In the next step one
replaces  the finite subset $E$ by an increasing sequence
 $(\Lambda_n)$ of  finite subsets of $\Gamma$, whose union is $\Gamma$ and
one proves that the sequence of operators $Op_h^{hyb , \Lambda_n}
(F)$ is a Cauchy sequence in ${\cal L}({\cal H} (\Gamma))$ (Theorem
5.4). To establish an upper bound on  the norms of the operators
$Op_h^{hyb , \Lambda_n} (F)$ or
 $Op_h^{hyb , \Lambda_n} (F)- Op_h^{hyb , \Lambda_m} (F)$ one adapts, up to some
 details,
the integration by parts method on which one proof of the
  Calder\'on-Vaillancourt Theorem is based. More precisely, one adapts the
proof due to  A. Unterberger [U2], which relies on coherent states.

\bigskip

The limit of the sequence of operators will be denoted by
 $Op_h^{weyl}(F)$ and can be considered as the Weyl operator associated with
$F$. Under the hypothesis  $H_2(M, \varepsilon)$, with $h>0$ and a
summable family  $(\varepsilon_j)_{(j\in \Gamma)}$ ,  its norm will
satisfy
$$ \Vert Op_h^{weyl}(F) \Vert _{{\cal L} ({\cal H}(\Gamma))} \leq M \prod
_{j\in \Gamma} (1 + 225 \pi K_2 \sqrt {h} \varepsilon _j) \ , \leqno
(1.9)$$
where  $K_2= \sup _{j\in \Gamma} \ \max (1 , \varepsilon_j^3)$.
Remark that, if the family  $(\varepsilon_j)_{(j\in \Gamma)}$ is
summable, then the infinite product converges. Under the hypothesis
$H_4(M, \varepsilon)$ of Definition 1.3, it is sufficient that
 $(\varepsilon_j^2)_{(j\in \Gamma)}$ be summable. If $0
< h \leq 1$, one has
$$ \Vert Op_h^{weyl}(F) \Vert _{{\cal L} ({\cal H}(\Gamma))} \leq M \prod
_{j\in \Gamma} (1 + 225 \pi K_4 h \varepsilon _j^2)\ ,  \leqno
(1.9')$$
where $K_4= \sup _{j\in \Gamma} \ \max (1 , \varepsilon_j^6)$.
Theorem 5.4 is the analog of Theorem 1.1 in an infinite dimensional
setting and its proof (see Sections 6, 7, 8) contains that of
Theorem 1.1.

\bigskip

One can wonder as to the relationship between this definition and
former definitions of the Weyl calculus, used for example by  B.
Lascar [LA1]
 and Kree-R\c aczka [K-R] or Albeverio-Daletskii [A-D].
Those definitions suppose in general that the symbol $F$ is the
Fourier
 transform of a finite measure on the Hilbert space  $\ell^2(\Gamma , \R)
\times \ell^2(\Gamma , \R)$. On the other hand, there is no
condition as to
 the existence and boundedness of partial derivatives.
We shall see in Section 9 that, if a symbol $F$  is the Fourier
 transform of a finite measure on  $\ell^2(\Gamma , \R)
\times \ell^2(\Gamma , \R)$ as well as verifies $H_2(M,
\varepsilon)$ (with $(\varepsilon_j)_{(j\in E)}$ a summable family),
then both definitions coincide (ours and, for example, B. Lascar's).
It would be interesting also to compare these definitions to
Khrennikov's [KH].

\bigskip

We now give an example of symbol satisfying our hypotheses.

\bigskip

{\it Example 1.4.} Let  $\Gamma = \Z^d$. For all  $j\in \Gamma$, set
$b_j = (1+|j|)^{\gamma}$, where $\gamma >0$. Let $B(\Gamma)$ be the
space of Definition  1.2, associated with this family. For all $X =
(x , \xi)$ in $B(\Gamma) \times B(\Gamma)$, set
$$ H(x , \xi ) = \sum _{j\in \Gamma} g_j^2 (x_j^2 + \xi_j^2) + \lambda
\sum _{|j-k|= 1} g_j g_k x_j x_k \ .  \leqno (1.10)$$
Here the norm  $|\cdot |$ is an arbitrary norm on  $\R^d$. If it is
the supremum norm  $\ell^{\infty}$, the function $H$ recalls a
lattice of harmonic oscillators with a coupling between nearest
neighbors. The constant $\lambda$ is such that the quadratic form
$H$ is positive definite and $(g_j)_{(j\in \Gamma)}$ is a family of
positive numbers. Set
$$ F(x , \xi) = e^{-H(x , \xi)} \hskip 2cm (x , \xi) \in B(\Gamma)
 \times B(\Gamma)\ .
\leqno (1.11) $$
If the family  $(g_j)_{(j\in \Gamma)}$ satisfies
$$ \sum _{j\in \Gamma} g_j^2 ( 1 + |j|)^{2 \gamma} < +\infty \ ,$$
then the function $F$ is continuous and bounded on $B(\Gamma)\times
B(\Gamma)$. Moreover, for every integer $m$, the function $F$
satisfies hypothesis
 $H_m(1 , \varepsilon)$, with $M=1$ and
$\varepsilon _j = C_m g_j$, where  $C_m$ is a constant depending on
$m$. Under these hypotheses, the condition $H_4(1 , \varepsilon)$
holds and the family
 $(\varepsilon_j^2)_{(j\in \Gamma)}$ is summable.
By the second version of Theorem 5.4, the Weyl operator associated
with the function  (1.11) will be bounded in
 ${\cal H}(\Gamma)= {\cal F}_s ( \ell^2(\Gamma , {\bf C}))$.

\bigskip

Now let us give an application of Theorem 1.1.

\bigskip

{\it Example 1.5.} With the notations of Example 1.4, let $(E_N)$ be
an increasing sequence of finite subsets of the lattice $\Gamma$,
whose union is $\Gamma$. Let $V$ be a real valued bounded function
in $C^{\infty}(\R)$,
 whose derivatives are all bounded. For every  integer $N$,
set
$$ H_N (x , \xi) = \sum  _{j\in E_N} \xi_j^2 + \sum _{(j , k)\in E_N
\times E_N \atop |j-k|\leq 1} V(x_j - x_k)\ , $$
$$ P_N (x , \xi) = e^{-{1\over |E_N|} H_N (x , \xi)}\ .$$
The function  $P_N$ satisfies condition  a) of Theorem 1.1, with
 $n = |E_N|$ and $M=1$. Condition b) is satisfied for the multi-indices
 $(\alpha , \beta)$ in $I_4( E_N)$, with
$\varepsilon_j = C_1 |E_N|^{-1/2}$, where  $C_1$ is a real constant.
By the second version of Theorem 1.1, the norm of $Op_1^{weyl}(P_N)$
satisfies :
$$ \Vert Op_1^{weyl}(P_N) \Vert _{{\cal L}(L^2 (\R^{E_N}))} \leq
\left ( 1 + {C_2 \over |E_N|} \right )^{|E_N|}\ ,$$
where $C_2>0$ is a constant. This norm has therefore an upper bound
which is independant of $N$.

\bigskip

Sections 2, 3 and 4 present the more or less classical notions,
which will be needed to state the main theorem: abstract Wiener
measures in Section 2, Fock spaces and Segal isomorphisms in Section
3, Segal Bargmann transform in an infinite dimensional situation in
Section 4. The account about the Segal Bargmann transform goes
further than needed in the rest of the article, but we thought it
was useful to clarify the connections between the different
definitions of this notion that can be found in the literature.
 In Section 5, we define the hybrid (Weyl- Anti Wick) operators
associated with the finite subsets of $\Gamma$ (Definition 5.1) and
we state Theorem  5.4, which can be considered as the main result.
In Section 6, we establish the formula linking two hybrid
quantizations associated with two finite subsets, one containing the
other. Sections 7 and 8 are devoted to the proof of the main result.
In Section 9, we compare our definition of the Weyl
pseudodifferential operator with the definition used in the articles
of Kree R\c aczka and B. Lascar.

\bigskip

{\bf 2. Measure spaces associated with the subsets of $\Gamma$. }

\bigskip

The Weyl operator associated with a symbol will be the limit of a
sequence of operators defined by a hybrid quantization, for which
the set $\Gamma$ has to be  split into complementary subsets  $E_1$
and $E_2$ playing different roles. This method compells us to
define, for each  subset $E$ of $\Gamma$, a configuration space
denoted by $B(E)$ and a phase space, which is naturally   $B(E)
\times B(E)$. For a finite $E$, $B(E) = \R^E$. Configuration and
phase spaces are equipped with  measures depending on a strictly
positive parameter $h$, respectively denoted by $\mu_{E,h}^K$ and
$\mu_{E,h}^{\Phi}$. For a finite $E$, the measure on $B(E)= \R^E$ is
$$ d\mu_{E,h}^K (u) = (\pi h)^{-|E|/2} \prod _{j\in E}
\left ( e^{-{1 \over h}u_j^2 } du_j \right ).\leqno (2.1)$$
The measure on the phase space  $\R^E\times \R^E$ is analogous to
 $\mu_{E,h}^K $ but, for technical reasons linked with the
coherent spaces, the variance is $2h$ instead of $h$. We have that

$$ d\mu_{E,h}^{\Phi} (x, \xi) = (2 \pi h)^{-|E|} \prod _{j\in E}
\left ( e^{-{1 \over 2 h}(x_j^2 + \xi_j^2) } dx_j d\xi_j \right )\ .
\leqno (2.2)$$

\bigskip
If $E$ is infinite, the definition of both measures will use the
theory of abstract Wiener spaces, [G2], [G3], Kuo [KU], Th. L\'evy
[LEV]). In order to apply the classical theorems allowing to define
measures
 analogous to  (2.1) and (2.2) on infinite dimensional spaces, we need
the configuration space $B(E)$ and its Banach space norm to satisfy
certain conditions. In the case when $E = \Gamma$, we shall find the
space $B(\Gamma)$ from Definition 1.2   again.

\bigskip

Before defining the convenient phase spaces $B(E)$, we  shall recall
the classical results needed. Let  $Z_{\R}$ be a separable Hilbert
space with norm $| \cdot |$. Let $B$ be a separable Banach space
with norm  $\Vert \cdot \Vert _B$, into which  $Z_{\R}$ is
continuously embedded as a dense subset. In this case,
  $B' \subset Z_{\R}$.
A{\it  cylinder set (or tame set)}
  of $B$ is a subset $X \subset B$ of the
 form
$$ X = \{ x\in B, \hskip 1cm \Big ( f_1(x) , ... , f_n(x) \Big ) \in
\Omega \}, \leqno (2.3)$$
where $n\geq 1$,  $f_1$, ... $f_n$ belong to  $B'$ and $\Omega$ is a
Borel set  of $\R^n$. A {\it cylinder (or tame)} function on $B$ is
a function  $f:B \rightarrow \R$ such that there exist elements
$f_1$, ... $f_n$ of $B'$ $(n\geq 1)$ and a function $g: \R^n
\rightarrow \R$ measurable for the Borel-$\sigma$ algebra  of $ \R^n
$ with which we can express $f(x)$ as
 $f(x) = g \left (f_1(x) , ... , f_n(x) \right )$
for all $x\in B$. Since the  $f_j$ can be considered as elements of
$Z_{\R}$, one can always assume that in (2.3), they form an
orthonormal system of $Z_{\R}$. In this case,
 if $X$ is a cylinder set defined by (2.3),  one sets
$$ \mu_{t, B} (X) = (2 \pi t)^{-n/2} \int _{\Omega} e^{-{|x|^2 \over
2t}} dx \leqno (2.4)$$
for all $t>0$.

This defines an additive mapping $\mu _{t ,B}$ on the set of all
cylinder sets of $B$.

 One can define the notion of  a cylinder set
in  $Z_{\R}$, using the same formula (2.3), but this time the  $f_j$
belong to  $Z_{\R}$ and $f_j(x)$ is the scalar product of $f_j$ and
$x$ in $Z_{\R}$. Similarly, the  formula (2.4) gives the
 measure of $X$
 (provided the $f_j$ are orthonormal). This defines an additive mapping
 on the set of all
cylinder sets of $Z_{\R}$. One shows (cf. Kuo [KU]) that, if
$Z_{\R}$ is infinite dimensional, the mapping  $\mu_{t, Z_{\R}}$
does not extend to a  $\sigma$-additive measure on the
$\sigma$-algebra generated by the cylinder sets of $Z_{\R}$. On the
other hand, the additive mapping $ \mu_{t, B}$ extends to a
$\sigma$-additive measure on the $\sigma$-algebra generated by the
cylinder sets of $B$, under hypotheses that we shall state now.

\bigskip

 A semi-norm   $N$ on $Z_{\R}$ is {\it   $\mu_{t, Z_{\R}}-$measurable}
if, for every number  $\varepsilon >0$, there exists a finite
 dimensional subpace $H_{\varepsilon}$ of $Z_{\R}$ such that, for
 every
finite
 dimensional subpace $V$ of $Z_{\R}$ orthogonal to $H_{\varepsilon}$,
the following inequality holds
$$ \mu _{t , Z_{\R}} \left ( \{ x \in Z_{\R} , \ \ \ \ N( P_V (x) )
> \varepsilon \} \right ) < \varepsilon\ , \leqno (2.5) $$
where $P_V$ is the orthogonal projection onto $V$ (so the set
appearing in (2.5) is  a cylinder set).

\bigskip

{\bf Theorem  2.1. (L. Gross)} {\it If $Z_{\R}$ is a real, separable
Hilbert space and $B$ a real separable Banach space into which
$Z_{\R}$ is continuously embedded as a dense subset, then the
$\sigma$-algebra generated by the cylinder sets is the Borel
$\sigma$-algebra of $B$. Moreover, if the norm of $B$, restricted to
$Z_{\R}$, is
 $\mu_{t, Z_{\R}}-$measurable, then the mapping  $\mu_{t , B}$ defined by
(2.4) on the cylinder sets extends to a uniquely determined measure
on the Borel $\sigma$-algebra of $B$. }

\bigskip

Both assertions are proved in Kuo [KU]. The first one ([KU] Theorem
4.2) does not require the $\mu_{t, Z_{\R}}-$ measurability of the
norm of $B$. The second one is due to L.Gross [G2], (see also
  [KU], Theorem 4.1 and  L. Gross [G3], Section 2, Theorem 1).
The following Proposition shows that in the infinite dimensional
case,
  $Z_{\R}$ is $\mu_{t, B}$-negligible whereas
 $\mu_{t, B}(B)=1$.

\bigskip

{\bf Proposition 2.2.} {\it With the notations and under the
hypotheses
 of  Theorem 2.1., if    $Z_{\R}$ is infinite dimensional,
it is contained in a  $\mu_{t, B}$-null Borel set. }

\bigskip

{\it Proof.} Let  $(e_n)_{(n\geq 0)}$ be an orthonormal basis of
$Z_{\R}$ whose vectors belong to $B'\subset Z_{\R}$. Let
$$ A = \{ x \in B , \hskip 1cm \sum _{n\geq 0}  e_n(x) ^2 < +\infty \} .$$
Clearly, $A$ is a Borel set of $B$ containing $Z_{\R}$. One defines
a sequence  $(\varphi_N)$ of cylinder functions by
$$ \varphi_N (x) = e^{ - \sum_{k=0}^N e_k(x)^2} \hskip 2cm x\in B$$
and denote by  $\varphi $ its pointwise limit.
 By the dominated convergence theorem,
$$ \int _B \varphi (x) d\mu_{t,B}(x) = \lim _{N \rightarrow \infty}
 \int _B  \varphi_N  (x) d\mu_{t,B}(x) = \lim _{N \rightarrow \infty}
(1+ 2t)^{-(N+1)/2} = 0.$$
Since $\varphi \geq 0$ and is stricty positive on $A$, it follows
that
 $\mu_{t,B}(A)=0$.

\bigskip

After recalling these results, we can apply them to prove the
 following theorem.

\bigskip

{\bf Theorem 2.3.} {\it Let  $(b_j)_{(j\in \Gamma)}$ be a family of
strictly positive real numbers, such that the family
 $R_j(b_j ,
\varepsilon)$ $(j\in \Gamma)$ of (1.6) is summable for every
$\varepsilon>0$. For every subset
 $E\subseteq \Gamma$, let
$B(E)$ be the space of the Definition 1.2. In this case, the space
 $Z_{\R} =
\ell ^2(E , \R)$ is densely embedded  in $B(E)$. Moreover, for all
 $t>0$, the restriction to $Z_{\R}$ of the norm of $B(E)$ is $\mu_{t,
Z_{\R}}-$measurable. }

\bigskip

The proof of the last point uses the following result (cf. L. Gross
[G3], Theorem 1, page 95).

\bigskip

{\bf Theorem 2.4.}  {\it Let  $\Vert \cdot \Vert _n$ be an
increasing sequence of tame semi-norms on
 $Z_{\R}$. Let  $t>0$.
If, for all $\varepsilon >0$,
$$ \lim _{n\rightarrow \infty}\mu_{t, Z_{\R}}\left ( \{  x \in H , \ \ \ \
\Vert x \Vert _n \leq \varepsilon \} \right ) > 0 \ ,\leqno (2.6)$$
then $\lim _{n\rightarrow \infty } \Vert x \Vert _n $ exists for all
 $x\in Z_{\R}$
and the limit defines a
 $\mu_{t,Z_{\R}}-$measurable semi-norm. }

\bigskip

{\it  Proof  of  Theorem  2.3 (last point). }
 One chooses an increasing sequence of finite subsets of $E$,
whose union is $E$. One defines  an increasing sequence $( \Vert
\cdot\Vert _N)$
 of tame semi-norms on  $H= \ell^2(E , \R)$  by
setting, for all
 $x
= (x_{\lambda})_{(\lambda \in E)}$
$$ \Vert x \Vert _N = \sup _{j\in E_N}{|x_{j}| \over b_j} \ . \leqno (2.7)$$
For all  $\varepsilon >0$, the set
$$ C_N = \{ x \in Z_{\R} = \ell^2(E , \R), \ \ \ \ \Vert x \Vert _N \leq
\varepsilon \}$$
is a cylinder set of $H$. It can be written as in (2.3) with the
Borel set
 $\Omega= \prod_{j\in E_N}[-\varepsilon b_j ,
\varepsilon b_j]$.  Its $\mu_{t, H}$- measure is therefore
$$ \mu_{t, H} (C_N) =  \prod _{j\in E_N} (2\pi t)^{-1/2}
\int _{-\varepsilon b_j}^{\varepsilon b_j} e^{-{x^2 \over 2t}} dx\ .
$$
The sequence $ \mu_{t, H} (C_N)$ decreases to a nonnegative limit.
One has
$$ \mu_{t, H} (C_N) =  \prod _{j\in E_N}\left [ 1 - 2 (2\pi t)^{-1/2}
\int _{\varepsilon b_j}^{+\infty } e^{-{x^2 \over 2t}} dx \right ] =
 \prod _{j\in E_N}\left [ 1 - 2 (2\pi )^{-1/2} R_j( b_j ,{ \varepsilon
 \over \sqrt {t}} ) \right ]\ , $$
where $R_j (., . )$ is defined by (1.6). Since the factors of this
product are $>0$, the limit is $>0$ provided the family  $R_j (b_j ,
{\varepsilon \over \sqrt{ t}} )$  is summable. As a consequence, the
sequence  $( \Vert \cdot\Vert _N)$ satisfies the condition of
Theorem  2.4 and  the limit above, which is the restriction of the
norm of $B$ to $Z_{\R}$, is  $\mu_{t,Z_{\R}}-$measurable.

\hfill \carre

\bigskip

We shall now give an example of a family $(b_j)_{(j\in \Gamma)}$
satisfying the condition of  Definition 1.2,
 in the case when $\Gamma$ is a lattice of  $\R^d$.

\bigskip

{\bf Proposition 2.5.} {\it Let  $\Gamma $ be a lattice of  $\R^d$
(with a norm
 $| \cdot |$) and let $\gamma >0$ be given. For all  $j\in \Gamma$, set
$b_j = (1+|j|)^{\gamma}$. Then for all $\varepsilon >0$, the family
$R_j ( b_j , \varepsilon )$ $(j\in \Gamma)$ defined by (1.6) is
summable.

}

\bigskip

{\it  Proof .} For all $j\in \Gamma$ one has :
$$2 (2\pi )^{-1/2}
\int _{\varepsilon (1+|j|)^{\gamma}}^{+\infty } e^{-{x^2 \over 2}}
dx  \leq  e^{-{\varepsilon^2 (1+|j|)^{2\gamma} \over 4}}(2\pi
t)^{-1/2}\int _{\R} e^{-{x^2 \over 4}} dx = \sqrt {2}
e^{-{\varepsilon^2 (1+|j|)^{2\gamma} \over 4}} .$$
Since  $\Gamma$ is a lattice of $\R^d$, we know that the family
$$ \sum _{j \in \Gamma }e^{-{\varepsilon^2 (1+|j|)^{2\gamma} \over 4}}$$
is summable, for every $\varepsilon >0$ and every dimension $d$.

\bigskip

According to Theorems 2.1 and 2.3, for every subset $E$ of $\Gamma$,
if $B(E)$ is the space of Definition 1.2, then the mapping  $\mu_{t
, B(E)}$ defined by (2.4) on the cylinder sets of $B(E)$ extends to
a uniquely determined measure on the Borel $\sigma$-algebra, still
denoted by   $\mu_{t , B(E)}$. In this paper, the measure on  the
configuration space is the measure  $\mu_{h/2 , B(E)}$ defined as
above with $t=h/2$ and an arbitrary $h>0$. From now on, it will be
denoted by  $\mu_{E ,h}^K$. The measure on the phase space is $\mu_{h
, B(E)}\otimes \mu_{h , B(E)}$, this time with $t=h$.   It will be
denoted by $\mu_{E ,h} ^{\Phi}$ in the rest of the paper and will be
used to define every integral on the phase space. If $E$ is finite,
then  $B(E)= \R^E$ and both measures coincide with those
 defined by
 (2.1) and (2.2).

\bigskip

One can now give an explicit expression of the integral of a
function $f$ in $L^1(B(E) , \mu^K_{E,h})$, at least when $f$ is a
cylinder or tame function. Suppose  that there exist a family $(z_1
,..., z_n)$ in $B(E)'$, orthonormal with respect to the scalar
product of $Z_{\R} = \ell^2 (E , \R)$ and a Borel measurable
function
 $F : \R^n\rightarrow {\bf C}$, such that
$$ f(x) = F( z_1(x), ... z_n(x)) \hskip 2cm x\in B(E). $$
One has then
$$ \int _{B(E)} f(x) d\mu_{E,h}^K(x) = (\pi h)^{-n/2} \int _{\R^n} F(u_1 , ... ,
u_n) e^{- {|u|^2 \over h}} d\lambda (u)\ ,$$
where $\lambda$ is the Lebesgue measure on  $\R^n$. In Section 3,
density arguments will allow us to extend this definition to the
case when $f$ is not a tame function but belongs to  $L^1(B(E) ,
\mu^K_{E,h})$.

\bigskip

Let us remark that, if  $E = E_1 \cup E_2$ with disjoint $E_1$ and
$E_2$, one has  $B(E)= B(E_1) \times B(E_2)$ and, for example for
the phase space

$$ \mu_{E,h}^{\Phi} = \mu _{E_1,h}^{\Phi}  \otimes \mu _{E_2,h}^{\Phi} .$$
In particular, if $E_1$ is finite, then $B(E) = \R^{E_1} \times
B(E_2)$ and
$$ d\mu_{E,h}^{\Phi} (X_{E_1}, X _{E_2} ) = (2 \pi h)^{-|E_1|}
e^{-{|X_{E_1}|^2 \over 2h}} d\lambda _{E_1} (X_{E_1})
 d\mu _{E_2,h}^{\Phi} ( X _{E_2} )\ ,  \leqno (2.8) $$
where $\lambda_{E_1}$ is the Lebesgue measure on $\R^{E_1}\times
\R^{E_1}$.

\bigskip

{\it Gaussian Vectors.     }

\bigskip

One knows (cf. Th. L\'evy [LEV]) that, if a real Hilbert space
$Z_{\bf R}$, a Banach space $B$ and a real number $t>0$ satisfy the
hypothesis of Theorem 2.2, then the complexified
 $Z_{\bf C}$ of $Z_{\R}$ is isomorphic to a subspace of $L^2(B,\mu_{t,B})$.
The elements of  $Z_{\R}$ are sent on elements of  $L^2(B , \mu_{t,
B})$ generally called {\it Gaussian random variables.}  We shall now
give the precise form of this isomorphism in the case when $Z_{\R} =
\ell^2(E , \R)$, for a subset $E$ of  $\Gamma$.

\bigskip

{\bf  Theorem  2.6.} {\it Let $E$ be an infinite subspace of
$\Gamma$. Let $a$ belong to
 $\ell^2 (E , {\bf C})$. For every finite subset $F$ of $E$, set
$$ \ell _{a , F} (x) = \sum _{j\in F} a_j x_j \ ,\hskip 2cm E_{a , F} (x) =
e^{\ell_{a , F}(x)}\ .$$
Let  $(F_n)_{(n\geq 0)}$  be an increasing sequence of finite
subspaces of $E$, whose union is $E$. The sequences of functions
$(\ell_{a , F_n})$ and
 $(E_{a , F_n})$ are Cauchy sequences in $L^2(B(E), \mu_{E,h}^K)$.
Their limits, respectively denoted by $\ell_a$ and $E_a$, are
independent of the sequence  $(F_n)$. The function $\ell_a$ belongs
to
 $L^p(B(E), \mu_{E,h}^K)$ $(1 \leq p < + \infty)$ as well and the
sequence  $(\ell_{a , F_n})$ converges to $\ell_a$ in $L^p(B(E),
\mu_{E,h}^K)$. One has
$$ \Vert \ell_a \Vert _{L^2(B(E))}^2  = {h\over 2} \Vert a \Vert _{\ell^2(E)}^2
  \hskip 1cm {\rm and }  \hskip 1cm \Vert E_a \Vert
_{L^2(B(E))}^2 = e^{h \Vert {\rm Re } (a) \Vert ^2 _{\ell ^2(E)} }.
\leqno (2.9)$$
One can write $E_a(x) = e^{\ell_a(x)}$ too. The system of all
functions
 $\ell_a (x)^p$ ($p$ integer  $\geq 0$, $a$ in $\ell^2(E ,
\R)$) is total in $L^2(B(E) , \mu ^K_{E,h})$. }

\bigskip

For the Hilbert space  $\ell^2(E , \R) \times \ell^2(E , \R) $ and
the phase space $B(E) \times B(E)$ one proceeds similarly, setting
$$ \ell _{(a , b)} (x , \xi) = \ell_a(x) + \ell_b (\xi).$$
Since the measure on $B(E) $ is $\mu_{E,h}^{\Phi}$, we need to alter
the computations (2.9). Writing $Z_{\bf C}$ instead of  $\ell^2(E ,
{\bf C})$, one gets, for all  $a$ in $Z_{\bf C} \times Z_{\bf C}$:
$$ \Vert \ell_a \Vert _{L^2(B(E)\times B(E), \mu_{E,h}^{\Phi})}^2  = h
 \Vert a \Vert _{Z_{\bf C}\times Z_{\bf C}}^2\ ,
  \hskip 2cm \Vert E_a \Vert
_{L^2(B(E)\times B(E), \mu_{E,h}^{\Phi})}^2 = e^{2h \Vert {\rm Re }
(a) \Vert ^2 _{Z_{\R} \times Z_{\R}} }\ . $$

\bigskip

 {\it Proof.} If  $n < m$, an explicit computation shows that
$$\Vert \ell_{a , F_n}  \Vert _{L^2(B(E))}^2  = {h\over 2}
\sum _{j\in F_n} |a_j|^2\ ,
 \hskip 2cm
\Vert \ell_{a , F_n} -  \ell_{a , F_m}  \Vert _{L^2(B(E))}^2 =
{h\over 2} \sum _{j\in F_m \setminus F_n} |a_j|^2.
$$
The assertions concerning the function $\ell_a$ easily follow from
this. One can see that
$$\Vert E_{a,F_n} \Vert _{L^2(B(E))}^2  = e^{ h\sum _{j\in F_n} |{\rm Re}(a_j)|^2} $$
and that, if $n <m$:
$$\Vert E_{a , F_n} - E_{a , F_m}  \Vert _{L^2(B(E))}^2  = ... $$
$$ ... = e^{ h\sum _{j\in F_n} |{\rm Re}(a_j)|^2}  \left [ e^{h \sum _{k\in
F_m  \setminus F_n} |{\rm Re}(a_k)|^2 }  - e^{{h\over 4} \sum _{k\in
F_m  \setminus F_n} a_k^2 } - e^{{h\over 4} \sum _{k\in F_m
\setminus F_n}\overline { a_k}^2 } + 1 \right ]. $$
Remarking that, for all $z\in {\bf C}$, one has $\left | e^z -1
\right | \leq |z| e^{ \sup ( {\rm Re}(z), 0)} $, one gets
$$\Vert E_{a,F_n} - E_{a,F_m}\Vert _{L^2(B(E))}^2  \leq 3h e^{h \Vert {\rm Re}(a)
\Vert _{\ell^2(E)}^2} \ \sum _{j\in F_m \setminus F_n} |a_j|^2 , $$
from which one can easily deduce the assertions concerning  $E_a$.
For the spaces  $L^p(B(E))$ it amounts to showing it for a real
number  $a$  and, via  H\"older, for  an even integer  $p \geq2$. In
this case one uses the inequality
 ${|u|^p \over p!}
\leq {1\over 2}(e^u + e^{-u}) -1$. One can pass from this to
$$ {1 \over p! } \int _{B(E)} |\ell_{a , F_m}(x) - \ell _{a ,
F_n}(x)|^p d\mu_{E,h}^K(x) \leq {h\over 4}e^{ {h\over 4}\Vert a\Vert
_{\ell^2(E)}^2} \sum _{j\in F_m \setminus F_n} |a_j|^2\ . $$
The results about the convergence in $L^p(B(E))$ are a consequence
of this.
 The result about the total system is classical
(cf. Janson [J]). It is based on the fact that, if a function $f$
belonging to  $L^2(B(E) , \mu^K_{E,h}) $ were orthogonal to all
$\ell^a (x)^p$ ($p$ integer $\geq 0$, $a$ in $\ell^2(E , \R)$), it
would be orthogonal to all functions  $x \rightarrow e^{i\ell _a
(x)}$. Accordingly, its Fourier transform with respect to the
Gaussian measure would be identically zero.

\bigskip

{\bf 3. Fock spaces and Segal isomorphisms.}

\bigskip

The definition of the hybrid Weyl-anti-Wick  quantization stated in
Section 5 involves the Hilbert space $L^2 (B(E)  , \mu_{E,h}^{K})$
(the configuration space) together with the Hilbert space $L^2 (B(E)
\times B(E) , \mu_{E,h}^{\Phi})$ (the phase space), for every subset
$E$ of $\Gamma$, where $B(E)$ is given by Definition 1.2 and where
the two measures are defined in Section 2.

\bigskip

Each of these two Hilbert spaces (the configuration and the phase
space) is isomorphic to some more abstract Hilbert space, namely,
the Fock space. The two corresponding Fock spaces will be used in
Section 4 in order to simplify its content.

\bigskip

{\it A. Symmetric Fock space over an Hilbert space.}

\bigskip

Even if  symmetric Fock spaces over  Hilbert spaces are standard
(c.f. [RE-SI], [SI1], [SI2], [J]), it may be useful to
recall their definitions here. Let $Z_{\bf C}$ be a separable complex Hilbert
space. For all $n\geq 2$, we denote by
$Z_{\bf C}^{\otimes n}$ the  $n$-fold tensor product
$Z_{\bf C}\otimes\cdots\otimes Z_{\bf C}$. The scalar product of
$X = u_1 \otimes ... \otimes u_n$ and $Y = v_1 \otimes ... \otimes
v_n$ belonging to  $Z_{\bf C}^{\otimes n}$ ($u_j \in Z$, $v_j \in Z$) is defined by:
$$ < X , Y > = < u_1 , v_1> ... <u_n , v_n> \ . \leqno (3.1)$$
The set  $Z_{\bf C}^{\odot n}$ is the subspace of $Z_{\bf
C}^{\otimes n}$ containing only symmetric elements. This space is
also often  designated by $S_n Z_{\bf C}^n$. The subspace $Z_{\bf
C}^{\odot n}$ of $Z_{\bf C}^{\otimes n}$ is associated with the
above scalar product when restricted to $Z_{\bf C}^{\odot n}$.

\bigskip

We agree that $Z_{\bf C}^{\odot 0} = {\bf C}$, and that $Z_{\bf
C}^{\odot 1}= Z$. The real number $1$ is considered as an element of
$Z_{\bf C}^{\odot 0}$, it is denoted by $\Omega$ and  called the
vacuum state. The algebraic direct sum of the   $Z_{\bf C}^{\odot
n}$ ($n\geq 0$) is associated with the scalar product satisfying
that the spaces $Z_{\bf C}^{\odot n}$ are $Z_{\bf C}^{\odot p}$ are
orthogonal when $n\not = p$. The completion of this direct sum under
the above scalar product is denoted by ${\cal F}_s(Z_{\bf C})$ and
is called the symmetric Fock space over $Z_{\bf C}$.

\bigskip

The symmetrized tensor product of  $n$ elements $u_1, \dots,
u_n$ in $Z_{\bf C}$ is defined by the following element of $Z_{\bf C}^{\odot n}$:
$$u_1 \odot\cdots  \odot u_n = {1 \over \sqrt {n!}} \sum _{\sigma \in
S_n} u_{\sigma (1)} \otimes \cdots \otimes u_{\sigma (n)}\ . \leqno (3.2)
$$
This notation is borrowed to  Th. Levy [LEV] and is different of those in
Reed-Simon [RE-SI]. The space $Z_{\bf C}^{\odot n}$ is generated by these symmetrized tensor products. The scalar product of two elements in $Z_{\bf C}^{\odot n}$ and $Z_{\bf C}^{\odot m}$  vanishes if $n\not = m$. When $m=n$,
$$ \Big < u_1 \odot\cdots \odot u_n  , v_1 \odot\cdots  \odot v_n  \Big
> = \sum _{\sigma \in S_n} \prod _{j=1}^n < u_j , v_{\sigma (j)}> 
\ .\leqno (3.3)$$
For each $X\in Z_{\bf C}$, the creation operator
$a^{\star}(X)$ acts in the algebraic direct sum of the
$Z_{\bf C}^{\odot n}$ and is defined by:
$$a^{\star}(X) ( u_1 \odot \cdots  \odot u_n ) = X \odot u_1 \odot u_2
\cdots \odot u_n  \leqno (3.4)$$
for any $u_1$, ... ,$u_n$ in $Z_{\bf C}$ and for each $n\geq 1$. It
maps the vacuum state to $a^{\star}(X) ( \Omega) = X$. The
annihilation operator $a(X)$ associated with $X\in Z_{\bf C}$ is
given by:
$$a(X) ( u_1 \odot\cdots \odot u_n )= \sum _{j= 1}^n  < u_j , X> u_1
\odot \cdots\odot \widehat {u}_j\odot \cdots \odot u_n \leqno (3.5)$$
for each $n\geq 2$.
The term $u_j$ is omitted in  the symmetrized product in the above right hand-side.
For $n= 1$, we have $a(X) ( u_1) = <u_1 , X > \Omega$, and if $n=0$ then
$a(X) (\Omega)=0$. Note that the mapping $X \rightarrow
a^{\star}(X)$ is ${\bf C}-$ linear, whereas the mapping $X
\rightarrow a(X)$ is antilinear.

\bigskip

For all $X$ in $Z_{\bf C}$, the Segal field is the unbounded operator $\Phi_S (X)$  defined by:
$$ \Phi_S ( X) (A) = {1 \over \sqrt {2}} ( a(X) + a^{\star}(X))
(A) \leqno (3.6)$$
for all $A$ in the algebraic sum of the $Z_{\bf C}^{\odot n}$.
It is essentially self-adjoint and its self-adjoint extension is also denoted by
$\Phi_S (X)$. The element $e^{ i
\Phi_S(X)} \Omega$ $(X \in Z_{\bf C})$ verifies:
$$ e^{ i \Phi_S(X)} \Omega = \sum _{(n\geq 0)}
 {i^n e^{-{|X|^2 \over 4}} \over 2^{n/2} n!}X \odot \cdots \odot X  \ .
\leqno (3.7)$$
(The $ e^{ i \Phi_S(X)}$ are often called coherent states of  ${\cal
F}_s (Z_{\bf C})$.)

\bigskip

{\bf Proposition 3.1.} {\it The set of $ e^{ i \Phi_S(X)}
\Omega $, $X\in Z_{\bf C}$, is complete in ${\cal
F}_s (Z_{\bf C})$. }

\bigskip

This result is probably very well-known. To prove it, we note that the set of subspaces $Z_{\bf C}^{\odot n}$ $(n\geq 0)$
is  complete  in ${\cal F}_s (Z_{\bf C})$ and by polarization (c.f. Janson [J], Theorem D.1), the set of elements $X \odot
\cdots \odot X$ $(X\in Z_{\bf C})$ is complete in
$Z_{\bf C}^{\odot n}$. Next, we deduce from (3.7) that (c.f. Rodnianski
Schlein [RO-SC]):
$$ {1\over 2\pi} \int _0^{2\pi}e^{-in\theta} e^{ i \Phi_S(e^{i\theta } X)} \Omega
d\theta ={i^n e^{-{|X|^2 \over 4}} \over 2^{n/2} n!}X \odot \cdots
\odot X$$ for each $X$ in $Z_{\bf C}$ and each $n\geq 0$.
Proposition 3.1 thus easily follows.

\bigskip

{\it Application.} Let  $E$ be a subset of $\Gamma$.  The symmetric
Fock  spaces over the Hilbert spaces $Z_{\bf C} = \ell^2 (E, {\bf
C})$ and $Z_{\bf C} = (\ell^2(E , {\bf C}))^2$ are  denoted by
${\cal H}(E)$ and ${\cal H}_{\Phi}(E)$ respectively. These two
spaces shall be our configuration and phase spaces in the following
sections. The corresponding vacuum states are denoted  by $\Omega_K
(E)$ and $\Omega _{\Phi}(E)$ respectively. The spaces ${\cal
H}^{fin}(E)$ and ${\cal H}_{\Phi}^{fin}(E)$ stands for the algebraic
direct sum of the subspaces  $Z_{\bf C}^{\odot n}$. The (Hilbertian)
completion of the tensor product of two complex Hilbert spaces
${\cal H}_1$ and ${\cal H}_2$ is denoted by  ${\cal H}_1 \otimes
{\cal H}_2$. Let us recall that
$$ {\cal F}_s (Z_1 \oplus Z_2) \simeq {\cal F}_s (Z_1)
 \otimes {\cal F}_s (Z_2)\ , $$
where  $Z_1$ and $Z_2$ are two complex Hilbert spaces.
In particular, when $E_1$ and $E_2$ are two disjoint subsets of  $\Gamma$,
the previous identity applied with $Z_j = \ell^2
(E_j , {\bf C})$ and with $Z_j =  \ell^2(E_j , {\bf C})^2$, ($1
\leq j \leq 2$) gives the following two identities:
$$ {\cal H} (E_1 \cup  E_2) \simeq {\cal H}(E_1) \otimes {\cal H}(E_2)\ ,
\hskip 2cm  {\cal H}_{\Phi} (E_1 \cup  E_2) \simeq {\cal H}_{\Phi}(E_1)
\otimes {\cal H}_{\Phi}(E_2)\ . \leqno (3.8) $$

\bigskip

{\it B. The Segal isomorphism.}

\bigskip

It is known that the Fock space ${\cal F}_s (Z_{\bf C})$ is
isomorphic to $L^2 (B , \mu _{t ,B})$ if the real Hilbert space
$Z_{\R}$ and the Banach space $B$ satisfy the hypotheses of Theorem
2.1 and if
 $Z_{\bf C}$ is the complexification of
$Z_{\R}$. A proof of this result appears in Janson [J] (Theorem
4.1). The starting point in [J] relies on the fact that $Z_{\bf C}$
is isomorphic to a subspace of $L^2 (B , \mu _{t ,B})$. This point is
also proved in our case in Theorem 2.6. Let us now recall the
construction of this isomorphism when $Z_{\R} = \ell^2(E , \R)$ with
$E$ being any subset of $\Gamma$. This isomorphism shall be denoted
by $J_{Eh}^K$.

\bigskip

It is sufficient to define $J_{E h}^{K}$ restricted to the subspaces
$Z_{\bf C} ^{\odot n}$ ($n\geq 0$).

\bigskip

For $n=0$,   we set
$$J_{Eh}^K (\lambda \Omega_{K}(E) )= \lambda \hskip 2cm \lambda \in {\bf C}
\ , \leqno (3.9)$$
where $\Omega_{K}(E)$ is the vacuum state of ${\cal H}(E)$.
The right hand-side is the constant function defined on  $B(E) $ equal to $\lambda$.

\bigskip

For  $n=1$, for any  $u$ in $Z_{\bf C}^{\odot 1}= Z_{\bf C} =
\ell^2(E , {\bf C})$, we set
$$J_{Eh}^{K}(u) (x) = \sqrt {2 \over h}  \ell_u(x) \leqno (3.10)$$
for a.e. $x$ in $B(E)$, where $\ell_u$ is the function belonging to
$L^2(B (E)) , \mu_{E,h}^K)$ defined in Theorem 2.6.

\bigskip

For $n\geq 2$, it suffices to define $J_{Eh}^{K}( u_1 \odot  \cdots
\odot u_n)$ for  $u_1$, \dots, $u_n$ in $Z_{\bf C}= \ell^2(E , {\bf
C})$. Let $P_n$ be the subspace of $L^2(B(E), \mu_{E,h}^K)$ spanned
by $1$ and by the functions $\ell _{v_1}\dots\ell_{v_k}$ ($ v_1$,
\dots ,$v_k$ in $ Z_{\bf C}$, $k \leq n$). Here, the product is the
multiplication product for functions on $B(E)$ and it belongs to
$L^2(B(E) , \mu_{E,h}^K)$. Let $\Pi_n$ be the orthogonal projection
in $\overline {P_n}$ on the  orthogonal to $\overline {P_{n-1}}$,
with the scalar product of $L^2 (B(E) , \mu_{E,h}^K)$. Then set
$$ J_{Eh}^{K} ( u_1 \odot  \cdots \odot u_n) =
\Pi _n ( J_{Eh}^{K}(u_1) ... J_{Eh}^{K}(u_n)) \leqno (3.11) $$
for all $u_1, \dots , u_n$ in
$Z_{\bf C}$.
The function $J _{Eh}^{K}( u_1 \odot  \cdots \odot u_n )$ is denoted by $ : u_1
\cdots u_n:$ by Th. Levy [LEV] (when $h=1$).

\bigskip

The above construction is identical to the one of  Janson [J]. It is
also noticed that, the mapping $J_{Eh}^K$ is extended by density as
an isometric isomorphism between ${\cal H} (E)$ and $L^2(B(E),
\mu_{E,h}^K)$, for any $h>0$, and each (finite or infinite) subset
$E$ of $\Gamma$.

\bigskip

Now setting  $Z_{\bf C}  = \ell^2(E , {\bf C})^2$, we naturally
proceed similarly for the construction of the mapping
$J_{Eh}^{\Phi}$ associated with the phase space. The only difference
with the construction associated with the configuration space being
that,  (3.10) becomes: for all $(u , v)$ in $ \ell^2(E , {\bf
C})^2$,
$$J_{Eh}^{\Phi}(u , v) (x, y) = \sqrt {1 \over h} \Big [ \ell_u(x)
+  \ell_v(y) \Big ] \leqno (3.12)$$
for a.e. $(x , y)$ in $B(E) \times B(E)$. The function $\ell_u$ is
defined in Theorem 2.6. Note that  $2h$  appears for the phase space
case whereas it is $h$ in the configuration space case. The
remaining part of the construction is identical to the one
concerning the configuration space.

\bigskip

{\it C. Segal isomorphism and Hilbertian bases.}

\bigskip

Let $E$ be a subset of $\Gamma$. We shall define the two Hilbertian
bases of  ${\cal H}(E)$ and $L^2(B(E) , \mu_{E,h}^K)$ (resp. of
${\cal H}_{\Phi}(E)$ and $L^2(B(E) \times B(E) , \mu_{E,h}^{\Phi})$
).  The Segal isomorphism realizes a bijection between the two bases
associated with the configuration space (resp. with the phase
space). We first fix some notations concerning multi-indices.

\bigskip

We call  a {\it multi-index } any map  $\alpha $ from $E\subseteq \Gamma$
into $\N $ such that $\alpha_j = 0$ except for a finite number of
indices $j$. The sum of $\alpha_j $ $(j \in E)$ is denoted by
$|\alpha| $. Let $S(\alpha )$ be the largest subset $S$ of $E$ such
that $\alpha_j \geq 1$ for all $j\in S$. It is necessarily finite.
We set:
$$ \alpha! = \prod _{j\in S(\alpha )} \alpha_j!\  . $$
Let $ (e_j)_{(j\in E)}$ be the canonical basis of  $Z_{\bf C} =
\ell^2(E, {\bf C})$. For every multi-index $\alpha $,
$e^{\alpha}\in{\cal H}(E)$ stands for the symmetrized product of
$|\alpha| $ factors, where each factor in the symmetrized product is
an $e_j$ ($j\in S(\alpha ))$,  and where each $e_j$ appears exactly
$\alpha_j$ times in the  product.

\bigskip

Hermite polynomials is the sequence of polynomials $H_n$ on $\R$, being an orthogonal
basis of $L^2(\R , \nu)$, where
$\nu$ is the measure $(2 \pi) ^{-1/2} e^{-{x^2\over 2 }} dx$, $H_n$ is
of degree $n$, and the coefficient of $x^n$ in $H_n$
 equals to 1. The $L^2(\R , \nu)$ norm of $H_n$ is $\sqrt {n!}$.
 For every multi-index $\alpha$, we set:
$$ P_{\alpha  , h }(x  ) = \prod_{j\in S(\alpha ) }
H_{\alpha _j } \left (  x_j  \sqrt {2 \over h} \right )\ . \leqno
(3.13) $$

\bigskip

{\bf Proposition 3.2.} {\it Let $c_{\alpha} = (\alpha ! )^{-1/2}$,
for every multi-index $\alpha$. The set of $c_{\alpha } e^{\alpha}$
is an Hilbertian basis of ${\cal H}(E)$. The set  of $c_{\alpha}
P_{\alpha  , h }$ is  an Hilbertian basis of $L^2( B(E) , \mu
_{E,h}^{K})$. The isomorphism $J_{Eh}^K$ satisfies:
$$ J_{Eh}^K (e^{\alpha}) = P_{\alpha h} \leqno (3.14) $$
for every multi-index $\alpha$.

}

\bigskip

{\it Proof.} The fact that the family of $c_{\alpha }
e^{\alpha}$ is orthonormal  follows from (3.3). The fact that the set of functions $c_{\alpha} P_{\alpha  , h }$ is orthonormal in $L^2( B(E)  , \mu _{E,h}^{K})$  comes from the above properties on Hermite polynomials.
Let us now show that this system is complete. We know that the set of functions on $B(E)$
$$ f(X) = \ell _{a}(X)^p \leqno (3.15) $$
is complete in $L^2( B(E) , \mu _{E,h}^{K})$, where $p\in\N$, $a$
belongs to $Z_{\bf C} = \ell^2(E , {\bf C})$, and $\ell_a$ is
defined in Theorem 2.6 (c.f. Janson [J] or see Theorem 2.6). Let $f$
be a function written as in (3.15), with $a$  in $Z_{\bf C}$. There
is a sequence $(a^{(\nu)})_{(\nu \geq 0)}$ of elements in $Z_{\bf
C}$ such that $a^{(\nu)}$ is the vanishing sequence except for a
finite number of indices, and such that the sequence
 $(a^{(\nu)})_{(\nu \geq 0)}$ converges to $a$ in $Z_{\bf
C}$ as $\nu $ tends to $+ \infty$. The following function:
$$f^{(\nu)}(X) = ( \ell _{a^{(\nu)}}(X)) ^p   $$
is a polynomial depending on a finite number of variables. In
Theorem 2.6, it is seen that:
$$ \lim _{\nu \rightarrow +\infty} \Vert \ell _{a} - \ell _{a} ^{(\nu)}\Vert
_{L^{2p} (B(E), \mu_{E,h}^K)} = 0\ .$$
For functions $f$ as in (3.15), we deduce that:
$$ \lim _{\nu \rightarrow +\infty} \Vert f - f^{(\nu)} \Vert
_{L^{2} (B(E), \mu_{E,h}^K )} = 0\ . $$
Consequently, the set of polynomial functions depending on a finite
number of variables is dense in the space $L^2( B(E)  , \mu
_{Eh}^{K})$. Since each polynomial functions depending on a finite
number of variables is a linear combination of the  $P_{\alpha , h
}$, then the set of $c_{\alpha}P_{\alpha , h }$ is an Hilbertian
basis of $L^2( B(E)  , \mu _{E,h}^{K})$. We can similarly show that
the  set of $c_{\alpha}e^{\alpha}$ is an Hilbertian basis of ${\cal
H} (E)$. Equality (3.14) follows from (3.11), (3.10) and the
following identity:
$$\left ( {2 \over h} \right )^{|\alpha|/2} \left (  \Pi_n ( f_{\alpha})
\right ) (x) = \prod _{j\in S(\alpha ) } H_{\alpha _j} \left ( x _j
\sqrt {2 \over h} \right ) \hskip 2cm f_{\alpha} (x) = \prod _{j\in
S(\alpha)} x_j ^{\alpha _j} $$
for every  multi-index $\alpha$, with $n = |\alpha|$.

 \hfill \carre

\bigskip

In order to consider the phase space, we set $u_j = (e_j, 0) $ and
$v_j =(0, e_j)$, for all $j\in \Gamma$. For every multi-index
$(\alpha , \beta)$, we define
$u^{\alpha}v^{\beta}$ in the phase Fock space ${\cal H}_{\Phi}
(E)$ as the symmetrized product of $|\alpha| +|\beta|$ factors,
each of them being either one of the $u_j$ ($j\in S(\alpha , \beta))$
or one of the  $v_j$ $(j\in S(\alpha , \beta))$, where each $u_j$ appears
 $\alpha_j$ times and each
$v_j$ appears $\beta _j$ times in the symmetrized product. Set:
$$ P_{\alpha, \beta , h }(x , \xi ) = \prod_{j\in S(\alpha , \beta) }
H_{\alpha _j } \left ( { x_j \over \sqrt {h}} \right ) H_{\beta _j }
\left ( { \xi_j \over \sqrt {h}} \right ) \leqno (3.16)$$
for every multi-index $(\alpha , \beta)$.
\bigskip

{\bf Proposition 3.3.} {\it The set of $c_{\alpha}c_{\beta}
u^{\alpha}v^{\beta}$ (with $c_{\alpha} = (\alpha !)^{-1/2}$) is an
Hilbertian basis of the Fock space   ${\cal H}_{\Phi}(E)$. The set
of functions $c_{\alpha} c_{\beta}   P_{\alpha \beta , h }$ is an
Hilbertian basis of $L^2( B(E) \times B(E) , \mu _{E,h}^{\Phi})$. The
Segal isomorphism $ J^{\Phi }_{Eh}$ verifies:
$$ J^{\Phi }_{Eh}(u^{\alpha}v^{\beta})  =
P_{\alpha \beta , h }  \leqno (3.17)$$
for every multi-index $(\alpha , \beta)$.

}

\bigskip

The proof is the same as the one of proposition 3.2.

\bigskip

{\it D. Covariance formulas.}

\bigskip

In view of  identities (3.8) we shall specify the two isomorphisms
$J_{Eh}^K$ and $J_{Eh}^{\Phi}$ when $E = E_1 \cup E_2$ with disjoint
$E_1$ and
 $E_2$. Let $X_E = (X_{E_1} , X_{E_2})$
be the running variable in $\R^E \times \R^E$. We have:
$$ \Big ( J_{E h}^{\Phi} (f_1 \otimes f_2 ) \Big ) (X_{E_1} ,
X_{E_2}) =  \Big ( J_{E_1 h}^{\Phi} (f_1) \Big ) (X_{E_1}) \
 \Big ( J_{E_2 h}^{\Phi} (f_2) \Big ) (X_{E_2}) \leqno (3.18) $$
for every $f_1$ in ${\cal H}_{\Phi}(E_1)$ and every $f_2$ in ${\cal
 H}_{\Phi}(E_2)$.

\bigskip

{\bf Theorem 3.4.} {\it For any (finite or infinite) subset $E$ of $\Gamma$, for every $f$ in ${\cal H} (E)$, for each $h>0$
and for all $a+ ib$  in $Z_{\bf C}=  \ell^2(E , {\bf C})$, ($a$
and $b$ being real numbers), we have:
$$ \Big ( J^K_{Eh} e^{ i \Phi_S (a + ib)} f \Big )(u)= e^{ -{1 \over 2 }|b|^2
+ {i\over 2}a.b + {i\over \sqrt {h} } \ell_{a+ib} (u) } \Big
(J^K_{Eh} f \Big ) ( u + \sqrt {h} b) \leqno (3.19)$$
for a.e. $u$ in $B(E)$.
Similarly, for every $F$ in ${\cal H}_{\Phi} (E)$, for each $h>0$
and for all $(a+ ib, a' + i b')$  in $Z_{\bf C}\times Z_{\bf C}$,
for a.e. $(x , \xi )$ in $B(E)\times B(E)$, we have:
$$ \Big ( J^{\Phi}_{Eh} e^{ i \Phi_S ( a + i b, a' + i b')} F \Big )(x, \xi )=
e^{\psi}  \Big ( J^{\Phi}_{Eh} F \Big )(x + \sqrt {2h} b, \xi -
\sqrt {2h} b')\ , \leqno (3.20)$$
$$ \psi = -{1 \over 2 }(|b|^2+ |b'|^2)
+ {i\over 2}(a.b + a'.b') + {i\over \sqrt {2h} }( \ell_{a+i b} (x) +
\ell_{a'+ i b'} (\xi )  )\ . \leqno (3.21)
$$

}

\bigskip

{\it Proof. First step:} we prove here
(3.19) for $f = \Omega _K(E)$. Let $X= a+i b$ be in $Z_{\bf
C}$. The isomorphism  $J^K_{Eh}$ is applied to  both sides of equality
(3.7). We use the definition (3.10)-(3.11) of this
isomorphism together with the following notations (c.f. Janson [J]): the element $ :
\ell _{u_1} \cdots\ell _{u_n} : $ ($: \ell_u ^n:$ when $u_1= \dots =u_n = u$) stands for the function $\Pi _n ( \ell _{u_1} \cdots \ell _{u_n}) $.
It is called Wick product. We obtain, with obvious notations:
$$ J^K_{Eh} \left ( e^{ i \Phi_S(X)} \Omega_K(E) \right ) =e^{-{|X|^2 \over 4}}
 \sum _{n\geq 0}  {i^n  \over h^{n/2} n!} :\ell_X^n : =
e^{-{|X|^2 \over 4}} : e^{{i\over \sqrt {h}} \ell_X}: \ . \leqno (3.22)
$$
According to Janson [J] (Theorem 3.33):
$$: e^{{i\over \sqrt {h}} \ell_X}: =  e^{{ i \over \sqrt {h}}\ell _X +
{1 \over 2h} E ( \ell_X^2)}\ ,  $$
where, with standard notations:
$$E ( \ell_X^2) = \int _{B(E)} ( \ell _a (u) + i \ell_b (u))^2 d\mu
_{Eh}^K (u)=  {h\over 2} (|a|^2 - |b|^2 + 2i a\cdot b )\ . $$
We then deduce (3.19) for $f = \Omega _K(E)$ since $J_{Eh}^K
(\Omega _K(E)) = 1$.

\smallskip

{\it Second step:} we now prove equality (3.19) when $f =
 e^{ i \Phi_S(Y)} \Omega_K(E) $, with $Y$ in $Z_{\bf C} =
\ell^2(E , {\bf C})$. From a standard formula on the
product of Weyl operators (c.f. Reed-Simon [RE-SI], Theorem
X.41, (X.65)), which is also a particular case of the Campbell
Hausdorff formula, we have:
$$ e^{ i \Phi_S(X)}  e^{ i \Phi_S(Y)} = e^{ {i\over 2} {\rm Im }( X
\cdot \overline Y) }e^{ i \Phi_S(X +Y)}\ . $$
Consequently:
$$ J_{Eh}^K \left (e^{ i \Phi_S(X)}  e^{ i \Phi_S(Y)}\Omega_K(E) \right )
 = e^{ {i\over 2} {\rm Im } (X
\cdot \overline Y) } \ J_{Eh}^K \left (  e^{ i \Phi_S(X
+Y)}\Omega_K(E) \right )\ . \leqno (3.23)$$
We apply the first step with $X$ replaced by $X+Y$ and with $X$ replaced by  $Y$. Combining these two formulas with (3.23) leads by direct computations to equality (3.19),
for $f =
 e^{ i \Phi_S(Y)} \Omega_K(E) $, with $Y$ in $Z_{\bf C} =
\ell^2(E , {\bf C})$. These elements
of ${\cal F}_s (Z_{\bf C} ) = {\cal H} (E)$ (coherent states)
form a total family in ${\cal H} (E)$ (proposition 3.1). Therefore, equality (3.19)
is valid for all  $f$ in ${\cal H}
(E)$. The proof of (3.20)(3.21) is a straightforward modification. Also note that the parameter $h$ becomes $2h$ when considering $J_{Eh}^{\Phi}$ instead of $J_{Eh}^K$.

\bigskip

{\bf 4. Segal Bargmann spaces and transforms.}

\bigskip

We shall use the coherent states, a standard family of functions of
$L^2(\R^E , \lambda _E)$, where $\lambda _E$ is the Lebesgue measure
on  $\R^E$, for every finite subset  $E$ of $\Gamma$. These functions
 may depend on the parameters $X =
(x , \xi)$ in $\R^E \times \R^E $ and $h>0$. These functions are here denoted
by $\Psi_{X , h}$ and are defined by:
$$\Psi^E_{X , h} (u) =  (\pi h)^{ -|E|/4}
e^{-{| u-x|^2 \over 2h}} e^{{i \over h} u .\xi - {i \over 2h} x.
\xi} \ , \hskip 2cm u\in \R^E\ ,   \leqno (4.1) $$
where the norm and the scalar product are those of $\R^E$. The exponent $E$
may be omitted from the notation. It is known that:
$$ < f , g> = (2\pi h)^{-|E|} \int _{(\R^E)^2} < f , \Psi^E_{X ,
h}> \ < \Psi^E_{X , h} , g> d\lambda_E(X) \leqno (4.2) $$
for all $f$ and $g$ in $ L^2(\R^E)$. In other word, one defines the
mapping  $\widetilde T_{Eh}$ from $L^2(\R^E)$ into $L^2(\R^E \times
\R^E)$ (where $\R^E$ is associated with the Lebesgue measure) by:
$$ (\widetilde T_{Eh} f) (x , \xi) = (2 \pi h)^{-|E| /2 } \int
_{\R^E} f(u) \overline {\Psi _{(x ,\xi ),h}(u)} d\lambda_E (u)\ . \leqno
(4.3) $$
Since our aim  will be to work in infinite dimension, we rather use the
gaussian measures   $\mu _{E,h}^K$ and
 $\mu _{E,h}^{\Phi}$ for the configuration and phase spaces, defined in
 Section 2. This leads us to define the following transform:
$$ ( T_{Eh} f) (x , \xi) =
\int _{\R^E} f(u) e^{ {1 \over h} u \cdot ( x - i \xi )  - {1 \over
4h} ( x - i \xi)^2 } d\mu _{E,h}^K (u) \ . \leqno (4.4)$$
This mapping is indeed a partial isometry from
 $L^2(\R^E,
\mu_{E,h}^K)$ into $L^2(\R^E \times \R^E, \mu_{E,h}^{\Phi})$, and  in finite dimension,
it is called  Segal Bargmann transform. One may see [FA][FO] for its properties in finite dimension. The range of this transform is the closed subspace of functions in $L^2(\R^E \times \R^E, \mu_{E,h}^{\Phi})$ which are
antiholomorphic once $\R^E \times \R^E$ is identified to ${\bf
C}^E$. This subspace is called  the Segal Bargmann space. We note that the integral transform in  (4.4) has been extended by  J. Sj\"ostrand [SJ] when the exponent in the right hand-side is a quadratic form or a more general function.

\bigskip

When the set $E$ is infinite, we shall also define a mapping, also
called the Segal Bargmann transform. It will be seen either as a
mapping from $L^2(B(E), \mu_{E,h}^K)$ into $L^2(B(E) \times B(E),
\mu_{E,h}^{\Phi})$, which extend the one in (4.4) or, as a mapping
$W_E$ from the configuration Fock space  ${\cal H}(E)$ into the
phase Fock space ${\cal H}_{\Phi}(E)$.  The two points of view are
equivalent according to the two Segal isomorphisms that are recalled
in Section 3. Some difficulties arise when defining  an analog to
the integral in (4.4) in infinite dimension. Therefore, we found
easier to define the Segal Bargmann transform  in the abstract Fock
spaces ${\cal H}(E)$ and ${\cal H}_{\Phi}(E)$.

\bigskip

{\it A. Definition  of the Segal Bargmann  transform.}

\bigskip

We shall first define this transform as a partial isometry $W_E$
from the configuration Fock space ${\cal H}(E)$ to the phase Fock
space ${\cal H}_{\Phi}(E)$ associated with $E$, for every (finite or
infinite) subset $E$.

\bigskip

We start by defining a partial isometry  $T$ from $Z_{\bf C} = \ell^2(E ,{\bf C})$
and taking values in  $Z_{\bf C} \times Z_{\bf C}$ by setting:
$$T ( u) = {1 \over \sqrt 2} ( u , -i u)\leqno (4.5)$$
for all $u\in Z_{\bf C}$
\bigskip

There is a canonical functor, usually denoted by $\Gamma$, which
associates to any continuous linear map $T$ from  a complex Hilbert
space $Z_1$ into a complex Hilbert space $Z_2$, with a norm smaller
than $1$, a map $\Gamma (T)$ from the Fock space ${\cal F}_s (Z_1)$
into ${\cal F}_s (Z_2)$, (see Derezinski G\'erard [D-G] (lemma 2.6)
or Reed-Simon [RE-SI], Section X.7, when $Z_1 = Z_2$). If $T$ is a
partial isometry from $Z_1$ into $Z_2$, $\Gamma (T)$ is a partial
isometry from ${\cal F}_s (Z_1)$ into ${\cal F}_s (Z_2)$. We shall
recall the definition of  $\Gamma (T)$ assuming (using the notations
in Section 3), that  $Z_1 = Z_{\bf C} = \ell^2(E ,{\bf C})$ and $Z_2
= Z_{\bf C} \times Z_{\bf C}$.

\bigskip

It is sufficient to define  $\Gamma (T)$ restricted to any subspace
$Z_{\bf C}^{\odot n}$ ($n\geq 0$). For $n=0$, we set $ \Gamma (T) (
\Omega _K(E)) = \Omega_{\Phi}(E)$. For $n=1$, we set $ \Gamma (T)
(u) = T(u)$, for every $u$ in $Z_{\bf C}$. When $n\geq 2$, let:
$$ \Gamma (T)( u_1 \odot \cdots \odot u_n)=T( u_1) \odot \cdots
\odot T( u_n) \leqno (4.6)$$
for all $u_1$,..., $u_n$ in $Z_{\bf C}$,
 ($n \geq 2$).
If $T$ has a norm smaller or equal than 1 then the mapping $\Gamma (T)$
defined above  is extended as a continuous linear mapping from
 ${\cal F}_s (Z_{\bf C}) = {\cal H}(E)$ into ${\cal F}_s (Z_{\bf C} \times Z_{\bf C})
 = {\cal H}_{\Phi}(E)$.

\bigskip

{\bf Definition 4.1.} {\it For every subset $E$ of $\Gamma$, we call
Segal Bargmann transform associated with $E$, the mapping $W_E$ from
${\cal H}(E)= {\cal F}_s ( Z_{\bf C}) $ ($Z_{\bf C}  = \ell^2(E
,{\bf C})$) into ${\cal H}_{\Phi}(E)= {\cal F}_s ( Z_{\bf C}\times
Z_{\bf C}) $, defined by  $W_E = \Gamma (T)$ where $T$ is the
mapping from $Z_{\bf C}$ into $Z_{\bf C} \times Z_{\bf C}$ defined
in (4.5).

}
 \bigskip

We note that:
$$ W_{E_1 \cup E_2} = W_{E_1} \otimes W_{E_2} \leqno (4.7)$$
when $E = E_1 \cup E_2$ with disjoint $E_1$ and $E_2$.

\bigskip

One obtains a partial isometry $\theta_{Eh}$ from $L^2(B(E) , \mu ^{K}_{Eh} )$ in $L^2(B(E)
 \times B(E) , \mu ^{\Phi}_{E,h} )$ defined by:
$$\theta_{Eh} = J_{Eh}^{\Phi} \circ W_E \circ \Big ( J_{Eh}^{K} \Big
)^{-1}  \leqno (4.8)$$
when composing $W_E$ defined above   with the two Segal isomorphisms.

\bigskip

The Lebesgue measure $\lambda _E$ becomes available again when $E$ is finite
and one may define the isomorphisms  $\widetilde J_{
Eh}^{K} $ and $\widetilde J_{ Eh}^{\Phi} $ between the configuration
Fock space ${\cal H}(E)$ (resp. phase Fock space ${\cal H}_{\Phi}(E)$) and
 $L^2(\R^E , \lambda_E)$ (resp. $L^2(\R^E
\times \R^E , \lambda_E\times \lambda _E)$). These isomorphisms are defined by:
$$ \widetilde J_{ Eh}^K(f) (u) =(\pi h)^{-|E|/4}  J_{Eh} ^ K(f) (u) \
e^{-{|u|^2 \over 2h}}   \hskip 2cm  u \in \R^E \leqno (4.9)$$
$$ \widetilde J_{ Eh}^{\Phi}(f) (x , \xi ) =(2 \pi h)^{-|E|/2}  J_{Eh} ^ {\Phi}(f) (x , \xi) \
e^{-{|x|^2+ |\xi|^2 \over 4h}}   \hskip 2cm  (x , \xi) \in \R^E
\times \R^E \ . \leqno (4.10)$$
Then, we can also define a mapping $\widetilde \theta _{Eh}$ from
$L^2(\R^E , \lambda_E)$ into $L^2(\R^E \times \R^E , \lambda_E\times
\lambda _E)$ by:
$$\widetilde \theta_{Eh} = \widetilde J_{Eh}^{\Phi} \circ W_E \circ
\Big ( \widetilde J_{Eh}^{K} \Big )^{-1} \ .  \leqno (4.11)$$
We shall verify (Theorem 4.3)  that, if $E$ is finite then  the two
mappings $T_{Eh}$ and $ \theta_{Eh}$ respectively defined in (4.4)
and (4.8) are equal. This holds true with the mappings $\widetilde
T_{Eh}$ and $\widetilde \theta_{Eh}$ defined (4.3) and (4.11).

\bigskip

{\it B. Segal Bargmann transform and  Hilbertian bases.}

\bigskip

Let $(e_j) _{(j\in E)}$ be the canonical basis of $Z = \ell^2(E, {\bf
C})$. Set
$$ w_j = {1 \over \sqrt 2} ( e_j , -i e_j )\ . \leqno (4.12) $$
For each multi-index $\alpha$, we define an element
$w^{\alpha}$ of the Fock space ${\cal H}_{\Phi}(E)$ as the symmetrized product of  $|\alpha|$ factors where each of the factor is a $w_j$ ($j$ in $S(\alpha)$) and where each factor $w_j$ $(j\in
S(\alpha))$ appears exactly  $\alpha_j$ times in the symmetrized product. We also set
$$ Q_{\alpha h  } (x , \xi) =  (2h)^{-|\alpha|/2}  \prod _{j\in S(\alpha)}
(x_j - i \xi_j)^{\alpha _j} \ . \leqno (4.13) $$

\bigskip

{\bf Proposition 4.2.} {\it Let $c_{\alpha} = (\alpha !)^{-1/2}$.
The set of elements $c_{\alpha} w^{\alpha}$ (where $S(\alpha )
\subseteq E$) is an orthonormal system of ${\cal H}_{\Phi}(E)$. The
set of functions $c_{\alpha} Q_{\alpha h}$ is an orthonormal system
in $L^2(B(E) \times B(E) , \mu_{E,h}^{\Phi})$. The $e^{\alpha}$ in
Section 3 satisfy:
$$ W_E  e^{\alpha} = w^{\alpha} \ . \leqno (4.14)$$
We have:
$$ J _{Eh}^{\Phi} w^{\alpha} = Q_{\alpha h} \ . \leqno (4.15)$$
}

\bigskip

{\it Proof.} Only  (4.15) needs to be proved.
We may write:
$$ w^{\alpha} = \sum _{\beta + \gamma = \alpha} { 2^{-|\alpha|/2} \alpha !
(-i)^{|\gamma| } \over \beta ! \gamma !} u^{\beta} v^{\gamma}\ ,  \leqno
(4.16)$$
where the  $ u^{\beta} v^{\gamma}$ are defined in Section 3. From
proposition 3.3, one has: $ J _{Eh}^{\Phi} (u^{\beta} v^{\gamma }) =
P_{\alpha \beta h}$ where $P_{\alpha \beta h}$ is defined in (3.16).
A basic formula on Hermite polynomials shows that:
$$Q_{\alpha h} =  \sum _{\beta + \gamma = \alpha} { 2^{-|\alpha|/2} \alpha !
(-i)^{|\gamma| } \over \beta ! \gamma !}P_{\alpha \beta h}\ . \leqno
(4.17)$$
The equality used, once written in one dimension,  is the following one and is probably standard:
$$ ( x- i \xi)^m = \sum _{p=0}^m C_m^p (-i)^{m-p} H_p(x)
H_{m-p}(\xi)\ . $$
Equality (4.15) then follows from (4.16),  proposition 3.3,
and (4.17).

\bigskip

{\it C. Integral form  (cylindrical case).}

\bigskip

{\bf Theorem 4.3.} {\it For every finite subset  $E$ of
$\Gamma$, the two mappings $T_{Eh}$ and $\theta _{Eh}$,
defined in (4.4) and (4.8), are equal. This is also true with  $\widetilde T_{Eh}$ and $\widetilde \theta _{Eh}$
defined in (4.3) and (4.11).  If $E$ is any (finite or infinite) subset of $\Gamma$ and if  $f\in L^2( B(E) ,
\mu^K_{E,h})$ depends only on the variables $u_j$ $(j\in S)$,
where $S$ is a finite subset of $E$, then:
$$ (\theta_{E h}f) (x , \xi) = \int _{\R^ S} f(u)e^{{1\over h} \varphi (
x , \xi, u)} d\mu _{E,h}^K (u)\ ,  \leqno (4.18)$$
where:
$$ \varphi (x ,\xi , u)= \sum _{j\in S} u_j( x_j - i \xi_j) -{1
\over 4} (x_j - i \xi_j)^2$$
For any finite $E$,  for every $X =(x , \xi)$  in $\R^E\times \R^E$,
let $\varphi_{X h}$ be the element ${\cal H}(E)$ defined by:
$$\varphi_{X h} = e^{ {i\over \sqrt {h}}\Phi_S ( \xi - i x) } \Omega_K(E)
\leqno (4.19) $$
with $\Phi_S$ given in (3.6). Let $ \Psi^E _{X h}$ be the coherent state defined in (4.1). Then, we have:
$$\widetilde J^{K}_{E h} \varphi_{X h} = \Psi_{x \xi h}^E  \ .
 \leqno (4.20) $$

}

\bigskip

{\it Proof.}  If  $S$ is finite then the mapping  $T _{Sh}$
defined in (4.4) is the Segal Bargmann transform in finite dimension
and its properties are well-known. It is a partial isometry from $L^2( \R^S , \mu _{S,h}^K)$ into $L^2 ( \R^S
\times \R^S , \mu _{S,h}^{\Phi})$ and verifies:
$$ T_{Sh} (P_{\alpha h} ) = Q_{\alpha h} \hskip 2cm S(\alpha) \subset S$$
where $P_{\alpha h}$ and $Q_{\alpha h}$ are defined in (3.13) and
(4.13). From propositions 3.2 and 4.2, the mapping $\theta
_{Sh}$ defined in (4.8) also satisfies:
$$ \theta _{Sh} (P_{\alpha h} ) = Q_{\alpha h} \hskip 2cm S(\alpha) \subset S$$
and it is also a partial isometry from $L^2( \R^S , \mu _{S,h}^K)$
into $L^2 ( \R^S \times \R^S , \mu _{S,h}^{\Phi})$. Thus, this two
mappings $T_{Eh}$ and $\theta _{Eh}$, defined in (4.4) and (4.8) are
equal. Consequently, the two mappings $\widetilde T_{Eh}$ and
$\widetilde \theta _{Eh}$ defined in (4.3) and (4.11) are also
equal. When $E$ is an arbitrary subset of  $\Gamma$ and if  $f$
belonging to $L^2( B(E) , \mu^K_{E,h})$ depends only on the variables
$u_j$ $(j\in S)$, where $S$ is a finite subset of $E$, then we see
that $\theta_{Eh}f$ depends only on the variables $(x_j , \xi_j)$
$(j\in S)$ and may be identified to the function $\theta _{Sh}f$,
and then to the function $T_{Sh}f$, which is the right hand-side of
(4.18). In order to prove (4.20), we apply Theorem 3.4 with $f$
replaced by $\Omega_K(E)$, $a$ by ${\xi \over \sqrt {h}}$ and $b$ by
$- {x \over \sqrt {h}}$, while taking into account that $J^K_{Eh}
\Omega _K (E) = 1$. We then immediately deduce (4.20) using the
definition (4.9) of the isomorphism $\widetilde J_{ Eh}^K$ and the
definition (4.1) of coherent states.

\bigskip

Theorem 4.3 provides another way to determine the mapping
 $\theta _{Eh}$ of (4.8) when  $E$ is infinite. For each $f$ in $L^2( \R^E , \mu _{E,h}^K)$, there is a sequence of functions
$f_n$ in the same space and a sequence of finite subspaces $S_n$ of $\Gamma$, such that $f_n $ depends only on the variables $u_j$ $(j\in S_n)$ and such that the sequence $(f_n)$
tends to $f$ in $L^2( \R^E , \mu _{E,h}^K)$. This fact follows from
proposition 3.2. The transforms $\theta _{Eh} (f_n)$ are determined by (4.18) with $S$ replaced by $S_n$. The sequence
$(\theta _{Eh} (f_n))$ is a Cauchy sequence in $L^2( \R^E , \mu
_{E,h}^K)$ and its limit is $\theta_{Eh}(f)$.

\bigskip

{\it D.  The Segal Bargmann space. }

\bigskip

We have to give a characterization of the range of the Fock space
 ${\cal H}(E)$ by
 $J_{Eh}^{\Phi} \circ W_E$ , that is to say, the range of the space
 $L^2( B(E) , \mu_{E,h}^K)$ by the operator $\theta_{Eh}$. It is a closed subspace of $L^2(B(E)
 \times B(E) , \mu ^{\Phi}_{E,h} )$. This point is well-known when $E$ is finite: the range is the subspace of $L^2(B(E)
 \times B(E) , \mu ^{\Phi}_{E,h} )$ constituted of antiholomorphic functions.
The article of B. Hall  [HA] remarks that the set of  (anti)holomorphic is not a closed subspace of
$L^2(B(E)
 \times B(E) , \mu ^{\Phi}_{E,h} )$. The Segal Bargmann space defined here has been introduced by  Driver-Hall [D-H].
It is defined in [D-H] as the $L^2( B(E)
\times B(E) , \mu _{E,h}^{\Phi})$ closure of the set of  antiholomorphic functions which are cylindrical, that it is to say, depending on a finite number of variables. In other words,
their definition is the property $i)$ below.

\bigskip

{\bf Theorem 4.4.} {\it Let $E$ be an infinite subset of
$\Gamma$ and $h>0$. Then, the following properties define equivalently the same closed subspace $SB(E , h)$ of $L^2(
B(E) \times B(E) , \mu _{E,h}^{\Phi})$:

\smallskip

i)  $SB(E , h)$ is the closure in $L^2( B(E) \times B(E) , \mu
_{E,h}^{\Phi})$ of the subspace of functions depending only on a finite number of variables, and which, additionnaly, are  antiholomorphic when identifying $B(E) \times B(E)$ with the complexification of  $B(E)$.

\smallskip

ii) $SB(E , h)$ is the closure in $L^2( B(E) \times B(E) , \mu
_{E,h}^{\Phi})$ of the subspace  spanned by the funtions
$Q_{\alpha , h}$ $(S(\alpha) \subseteq E)$.

\smallskip

iii) $SB(E , h)$ is the range of the configuration Fock space ${\cal H}(E)$ by the mapping $J_{Eh}^{\Phi}$.

\smallskip

iv) $SB(E , h)$  is the range  of the space
 $L^2(B(E) , \mu_{E,h}^K)$ by the mapping $\theta _{Eh}$ of (4.8).

}

\bigskip

{\it Proof.} In order to derive that the space defined in $i)$
is included into the one defined in $ii)$, we consider a function $f$, depending on a finite number of variables and being antiholomorphic. Taking  Proposition 3.3 into account, we may write:
$$ f = \sum _{S(\alpha , \beta) \subseteq E} f_{\alpha \beta}
c_{\alpha } c_{\beta} P_{\alpha \beta h}\ , \hskip 2cm
 \sum _{S(\alpha , \beta) \subseteq E} |f_{\alpha \beta}|^2 = \Vert
 f  \Vert ^2\ . $$
Let $\Pi$ be the  orthogonal projection on the subspace $SB(E , h)$
defined by  $i)$. Since $f = \Pi f$, we have:
$$ f = \lim _{N \rightarrow + \infty} f_N\ ,  \hskip 2cm
f_N = \sum _{S(\alpha , \beta) \subseteq E \atop |\alpha |+ |\beta|
\leq N}  f_{\alpha \beta} c_{\alpha } c_{\beta} \Pi P_{\alpha \beta
h}\ . 
$$
We see that  $\Pi P_{\alpha \beta h}$ is a linear combination
of the $Q_{\gamma , h}$ such that $|\gamma |\leq |\alpha|+|\beta|$ and
$S(\gamma ) \subseteq S(\alpha , \beta)$. Consequently, $f$ is in the space defined in  $ii)$. Thus, the space defined in  $i)$ is included in the space defined in  $ii)$. These two spaces are then equal.
We remark that the spaces defined in  $ii)$ and $iii)$ are equal from
proposition 4.2 using that the set of $c_{\alpha}e^{\alpha}$ is an Hilbertian basis of ${\cal H}(E)$ and since $J_{Eh}^{\Phi}
\circ W_E$ is a partial isometry. The equality between the spaces
defined in $iii)$ and  $iv)$ comes from the fact that  $J_{Eh}^{K}$ is an isomorphism from  ${\cal H}(E)$ to $L^2(B(E) , \mu_{E,h}^K)$.

\bigskip


{\bf Proposition 4.5. } {\it Let  $F$ and $G$ be in $SB (E , h)$.
Then, for every $a$ and $b$  in $Z_{\bf R} = \ell^2(E ,
\R)$, we have:
$$ \int _{B(E) \times B(E)} e^{ - {1 \over 2h} \ell_{a-ib} (x + i
\xi)} F( x +a , \xi +b) \overline {G(x , \xi )} d\mu_{E,h}^{\Phi}(x ,
\xi) = \int _{B(E) \times B(E)} F(X) \overline {G(X )}
d\mu_{E,h}^{\Phi}(X)\ .  \leqno (4.21)$$

}

\bigskip

{\it Proof.} We define the operator $T$ in
$L^2(B(E) \times B(E), \mu_{E,h}^{\Phi})$ by
$$ (T_{ab} \varphi  ) (x , \xi) =  e^{ - {1 \over 2h} \ell_{a-ib} (x + i
\xi)} \varphi ( x +a , \xi +b) \hskip 2cm \varphi \in L^2(B(E)
\times B(E), \mu_{E,h}^{\Phi})\ . $$
This operator is bounded in $L^2(B(E) \times B(E),
\mu_{E,h}^{\Phi})$ with a norm smaller than $e^{|Y|^2 \over 4h}$. Since the functions
$F$ and $G$ are in $SB(E , h)$, there are two sequences  $(F_n)$
and $(G_n)$ of cylindrical (depending on a finite number of variables) antiholomorphic functions  converging to
$F$ and $G$ in $L^2(B(E) \times B(E), \mu_{E,h}^{\Phi})$. For every  $n$, we have:
$$ < (T_{a b} - I) F_n , G_n > = 0\ .$$
Indeed, the proof of the proposition is elementary in finite dimension. According to the continuity of $T_{a b}$, we then deduce that $ < (T_{a b} - I) F , G > = 0$, which is Proposition 4.5.

\bigskip

The next issue is now to give an explicit integral expression of the
Bargmann transform in infinite dimension, which possibly extend
(4.4). We may, according to Section 4.C, approximate any arbitrary
function in $L^2(B(E) , \mu_{E,h}^K)$ by a a sequence $(f_n)$ of
cylindrical functions, on which we may apply the transform (4.18), which is
 the cylindrical analog of (4.4). There is also a transform
analog to (4.4) in  [K-R], but leaving the cylindrical framework.
This transform is used in  Lascar [LA1]. We shall see in Theorem
4.10 that in some sense, this integral transform is the transform
$\theta_{Eh}f$ defined in (4.8), once restricted to $Z_{\R} \times
Z_{\R}$. If  $E$ is infinite then $Z_{\R} \times Z_{\R}$ is of
measure zero in $B(E)\times B(E)$ (Proposition 2.2), and this notion
of restriction should be first clarified.

\bigskip

{\bf Theorem 4.6.} {\it For every subset $E$ of
$\Gamma$ and for each $X = (x , \xi)$ in $Z_{\R} \times Z_{\R}$,
the mapping $F \rightarrow F (X)$, defined on the space spanned by the $Q_{\alpha h}$ $(S(\alpha) \subseteq E)$,  can be extended in an unique way to a continuous linear form on the space  $SB(E , h)$. This extension is denoted by $\rho_X $. For every $F$ in $SB(E ,
h)$ and for any $X = (x , \xi)$ in $Z_{\R} \times Z_{\R}$, we have:
$$ \rho_X (F) =\int _{B(E) \times B(E)}e^{{1 \over 2h}
\ell _{x - i \xi } (y + i \eta) }\
 F(Y) d\mu ^{\Phi}_{E ,h} (Y)\ ,  \leqno (4.22)$$
where the above exponential is the function defined in Theorem 2.6.
One also has:
$$ \rho_X (F) =\int _{B(E) \times B(E)}e^{-{1 \over 2h}
\ell _{x + i \xi} (y - i \eta) }\
 F(X+ Y) d\mu ^{\Phi}_{E ,h} (Y)\ .  \leqno (4.23)$$
For each $F$ in $SB(E , h)$, the function $X\rightarrow \rho_X F$
is continuous on $Z_{\R} \times Z_{\R}$ and
G\^ateaux antiholomorphic when identifying $Z_{\R} \times Z_{\R}$ with $Z_{\bf C}$.

}

\bigskip

The integral operator  $\rho_X$ is often called {\it reproducing kernel}. However, this terminology seems appropriate only when $E$ is
finite, since in that case $Z_{\R} = B(E)= \R^E$.

\bigskip

{\it Proof of Theorem 4.6.}  For each $X$ in $Z_{\R} \times Z_{\R}
$, set:
$$ E_X (y , \eta) = e^{{1 \over 2h}\ell _{y - i \eta} (x + i \xi ) }\ . $$
Then, the mapping $X \rightarrow E_X$ is continuous from $Z_{\R}
\times Z_{\R}$ into $L^2(B(E) \times B(E) ,\mu
^{\Phi}_{E,h} )$ and we have:
$$ \Vert E_X \Vert  = e^{h\Vert {\rm Re} \ X \Vert _{\ell^2(E)}^2}\ . $$
Consequently, the integral in (4.22) properly defines a continuous linear form $\rho_X$ on $L^2(B(E) \times B(E) ,\mu
^{\Phi}_{E,h} )$.  The preceding remarks imply that the mapping
$X \rightarrow \rho_X(F)$ is continuous on $Z_{\R} \times Z_{\R} $ and G\^ateaux antiholomorphic, for each
$F$ in $L^2(B(E) \times B(E) ,\mu ^{\Phi}_{E,h} )$. Moreover, $\rho_X(F)$ can be expressed as in (4.23). If $F$ is a linear combination of the  $Q_{\alpha h}$ then there is a finite subset
$S$ of $E$ such that  $F$ depends only on the variables $x_j$ and
$\xi_j$ $(j\in S)$. In this situation, this also holds in the integral (4.22) which may be written as:
$$ \rho_X (F) = \rho_{X_S} (F) =(2 \pi h)^{-|S|} \int _{\R^S \times \R^S}e^{{1 \over 2h}
(x_S - i \xi_S ) \cdot  (y_S + i \eta_S) }\
 F(Y_S) e^{ - {1\over 2h} |Y_S|^2} \ d\lambda_S(Y_S)  \leqno (4.24)$$
for all $X$ in
$Z_{\R}\times Z_{\R}$ and where $\lambda_S$ is the Lebesgue measure on $\R^S \times \R^S$.
This integral makes sense since
 $S$ is finite. Since $F$
is identified to an antiholomorphic function on ${\bf C}^E$, then
the reproducing kernels theory in finite dimension shows that $
\rho_X (F) = F(X) $, for every function $F$ written as a finite
linear combination of the $Q_{\alpha h}$ and for each  $X$ in
$Z_{\R}\times Z_{\R}$. Therefore, $\rho_X$ defined by (4.22) or by
(4.23) is the unique extension by continuity of the mapping $F
\rightarrow F (X)$ defined by finite linear combinations of  the
$Q_{\alpha h}$.

\hfill \carre

\bigskip

In view of (4.22), we may give a fifth characterization of the space
$SB(E , h)$ equivalent to those in Theorem 4.4. For each finite
subset $S \subset E$, one may define a cylindrical reproducing
kernel
 $\rho^S$ in the following way. Set $T = E \setminus S$ and denote by $(X_S,
X_T)$ the variable in $B(E) \times B(E)$ with $X_S$ in $\R^S
\times \R^S$ and $X_T$ in $B(T) \times B(T)$. The operator
$\rho^S$ is defined by:
$$ \rho ^S F (X_S , X_T) =(2 \pi h)^{-|S|} \int _{\R^S \times \R^S}e^{{1 \over 2h}
(x_S - i \xi_S ) \cdot  (y_S + i \eta_S) }\
 F(Y_S, X_T ) e^{ - {1\over 2h} |Y_S|^2} \ d\lambda_S(Y_S)  $$
for all $F$ in $L^2(B(E) \times B(E),
\mu_{E,h}^{\Phi})$.
This operator is bounded in $L^2(B(E) \times B(E),
\mu_{E,h}^{\Phi})$.

\bigskip

{\bf Proposition 4.7.} {\it Let $E$ be a finite subset of $\Gamma$
and let $h>0$. The space $SB(E , h)$  in Theorem 4.4 is also
characterized by the following property:

\smallskip

v) $SB(E , h)$ is the set of all $F$ in $L^2(B(E) \times B(E),
\mu_{E,h}^{\Phi})$ satisfying $\rho^S F = F$ for every finite subset
$S$ of $E$.

}

\bigskip

{\it Proof.} Let $F$ be in the space defined by $v)$. For any
integer $N$ and every finite subset $S$ of $E$, let $\Pi_{NS}$ be
the orthogonal projection operator in $L^2(B(E) \times B(E),
\mu_{E,h}^{\Phi})$, on the subspace spanned by the $P_{\alpha \beta
h}$ with $|\alpha|+|\beta|\leq N$ and $S(\alpha , \beta) \subseteq
S$. For every $\varepsilon >0$, there are $N$ and a (finite) $S$
 satisfying $\Vert F - \Pi _{NS}F \Vert < \varepsilon $. We remark that $\rho ^S$ and $\Pi_{NS}$ commute. Consequently, $\Pi_{NS}F$
is a stable space by $\rho^S$. Thus, since $S$ is finite then
$\Pi_{NS}$ is a Bergman anti-projection and its range consists of
antiholomorphic functions. Thus, every function $F$ in the space
defined by $v)$ is limit of a sequence of antiholomorphic functions
depending only on a finite number of variables, and is therefore in
the space defined by $i)$ in Theorem 4.4.

\bigskip

{\it E. Covariance formulas.}

\bigskip

 {\bf Proposition 4.8.} {\it For any $f$ in ${\cal
H}(E)$, for any $Y = y + i \eta $ in $Z_{\bf C} = \ell^2(E ,{\bf
C} )$ and for a.e. $(x , \xi) $ in $B(E) \times B(E)$, we have:
$$\left ( J^{\Phi} _{Eh} W_E e^{{i\over \sqrt {h}} \Phi_S (iY)} f \right ) (x , \xi ) =
e^{-{1\over 4h}|Y|^2 -{1\over 2h} \ell _{y + i\eta}   ( x - i \xi) }
\left ( J^{\Phi} _{Eh} W_E f \right ) (x + y , \xi + \eta )\ . \leqno
(4.25)$$
}

\bigskip

{\it Proof.} From the definitions of the mapping
$W_E$ and of the Segal field $\Phi_S$:
$$ W_E \Phi_S (Y) f = {1\over \sqrt 2} \Phi_S ( Y ,  -i Y) W_E f \leqno (4.26)$$
for every $Y$ in $Z_{\bf C}$ and every $f$ in ${\cal H}^{fin}(E)$.
In the above right hand-side (resp. left hand-side), the Segal field
acts in ${\cal H}_{\Phi} (E)$, (resp. in  ${\cal H} (E)$). Taking
the exponential, we deduce that:
$$ W_E e^{ i \Phi_S (Y)}  f = e^{ i \Phi_S ({Y \over \sqrt 2} , - i {Y \over \sqrt 2} )}
W_E f   \leqno (4.27) $$
for all $f$ in
${\cal H}(E)$ and any $Y$ in $Z_{\bf C}$. Consequently,
$$\left ( J^{\Phi} _{Eh} W_E e^{{i\over \sqrt {h}} \Phi_S (iY)} f \right ) (x , \xi )
= \left ( J^{\Phi} _{Eh} e^{{i\over \sqrt {2h}} \Phi_S (iY, Y)}W_E f
\right ) (x , \xi)\ . $$
Applying Theorem 3.4    (point 2) with $F = W_Ef$ and $b = a' =
{y\over \sqrt {2h}}$, $b' = -a = { \eta \over \sqrt {2h}}$, we
obtain (4.25).

\bigskip

{\bf Theorem 4.9.} {\it For every subset $E$ of $\Gamma$, for every
$f$ and $g$ in ${\cal H}(E)$, for each $Y = y + i  \eta$ in $Z_{\bf
C} = \ell^2(E , {\bf C})$:
$$ < e^{-{i\over \sqrt {h} } \Phi_S (y + i \eta) } f , g > _{{\cal
H}(E)} = ... \leqno (4.28)$$
$$ ... = e^{{|Y|^2 \over 4h}}  \int _{B(E) \times B(E)}
e^{-{i\over h} \big (  \ell_y(x) + \ell_{\eta}(\xi) \big ) } \Big (
J_{Eh}^{\Phi} W_E f \Big )(x , \xi)  \overline { \Big (
J_{Eh}^{\Phi} W_E g \Big )(x , \xi) } d\mu _{E,h}^{\Phi} (x , \xi)\ .$$
}

\bigskip

{\it Proof.} Let $A$ be the left hand-side of (4.28). Applying the isometries to the two functions in the scalar product, we see that:
$$ A =\int _{B(E) \times B(E)}\Big (
J_{Eh}^{\Phi} W_E e^{-{i\over \sqrt {h} } \Phi_S (y + i \eta) } f
\Big )(x , \xi)  \overline { \Big ( J_{Eh}^{\Phi} W_E g \Big )(x ,
\xi) } d\mu _{E,h}^{\Phi} (x , \xi)\ . $$
From Proposition 4.8, $A$ satisfies:
$$ A =\int _{B(E) \times B(E)} \Phi_{Y} (x , \xi)  \overline { \Big ( J_{Eh}^{\Phi} W_E g \Big )(x ,
\xi) } d\mu _{E,h}^{\Phi} (x , \xi)\ , $$
$$\Phi_{Y} (x , \xi)  = e^{-{1 \over 4h}|Y|^2 - {i\over 2h} (\ell_y + i \ell _{\eta})
(x - i \xi)}\Big ( J_{Eh}^{\Phi} W_E f \Big )(x - \eta  , \xi + y )\ . 
$$
The function $\Phi$ is in $SB(E , h)$. The function
$G = J_{Eh}^{\Phi} W_E g$ is also  in $SB(E , h)$. Applying Proposition 4.5 to the functions $\Phi$ and $G$ with $a =\eta  $ and $b = -y $, we learn:
$$A =\int _{B(E) \times B(E)}
e^{ - {1 \over 2h} \ell_{\eta + iy} (x + i \xi)} \Phi_{Y} (x + \eta
, \xi - y )  \overline { \Big ( J_{Eh}^{\Phi} W_E g \Big )(x , \xi)
} d\mu _{E,h}^{\Phi} (x , \xi)\ . $$
A direct computation then gives (4.28).

\bigskip

{\it F. Connection to the definition of Kree R\c aczka. }

\bigskip

In Kree R\c aczka [K-R],  the Segal Bargmann transform is defined as
the mapping from  $L^2 ( B(E) , \mu _{E,h}^K)$ taking values in the
space of continuous functions being G\^ateaux antiholomorphic on
$Z_{\bf C} = \ell^2(E , {\bf C})$. This integral transform is
defined by the right hand-side below in equality (4.29). Identifying
$Z_{\bf C}$ and $Z_{\R} \times Z_{\R}$, we shall see that the
transform of  $f$  defined in [K-R] is actually also the function $X
\rightarrow \rho_X ( \theta _{Eh} (f))$, where $\theta _{Eh}$ is
defined in (4.8) and $\rho_X$ is defined in Theorem 4.6. We shall
also remark that any element $F$ in the Segal Bargmann space $SB(E ,
h)$ is uniquely determined by its "restriction" to $Z_{\R} \times
Z_{\R}$, that is to say, by the function $X \rightarrow \rho_X(F)$
defined on
 $Z_{\R} \times Z_{\R}$. In other words, the transform of $f$
defined in Kree R\c aczka [K-R] uniquely determines $\theta _{Eh} (f)$. The integral in (4.29) is the finite dimensional analog of  (4.4), but it makes sense, in general, only if $X = (x , \xi)$ is in $Z_{\R} \times Z_{\R}$ instead of
$B(E) \times B(E)$.

\bigskip

{\bf Theorem 4.10. }  {\it For every $f$ in $L^2(B(E) , \mu
_{E,h}^K)$ and for any $X = (x , \xi)$ in $Z_{\R} \times Z_{\R}$,
 $(Z_{\R} = \ell^2(E , \R))$, we have:
$$\rho _X \Big (\theta_{Eh} f \Big )  = \int _{B(E)} f(u) e^{ {1 \over h}
\ell _{ x - i \xi } (u) - {1 \over 4h} ( x - i \xi)^2 } d\mu _{E,h}^K
(u)\ .  \leqno (4.29)$$

}

\bigskip

The proof uses the proposition below.

\bigskip

{\bf Proposition  4.11.    } {\it  For any subset $E$, for every $f$ in ${\cal H}(E)$, for each $Y = (y , \eta)$
in $ \ell^2(E , \R)^2$ identified to $Y = y + i
\eta$ in $\ell^2(E , {\bf C})$, we have:
$$e^{-{1\over 4h} |Y|^2} \  \rho_Y (  J_{ Eh}^{\Phi}  W _{E}
f   )   = < e^{ {i\over \sqrt {h}} \Phi_S (iY) } f , \Omega_K(E) >\ . 
\leqno (4.30)$$

 }

 \bigskip

 {\it Proof.}    Fix $f$ in ${\cal H} (E)$.
 The function $F = J_{ Eh}^{\Phi}  W _{E} f$ belongs to $SB(E , h)$.
Take $Y$ in $ \ell^2(E , \R)^2$.
According to the second expression (4.23) of $ \rho_Y $, we note:
$$ \rho_Y \left (J_{ Eh}^{\Phi}  W _{E} f \right ) =\int _{B(E) \times B(E)}e^{-{1 \over 2h}
(\ell _y + i \ell_{\eta})(x - i \xi) }\
 \left ( J_{ Eh}^{\Phi}  W _{E} f \right ) (X+ Y) d\mu ^{\Phi}_{E,h} (X)
\ .  \leqno (4.31)$$
In view of Proposition 4.8, we obtain:
$$ e^{-{1 \over 2h}(\ell _y + i \ell_{\eta})(x - i \xi) }\
 \left ( J_{ Eh}^{\Phi}  W _{E} f \right ) (X+ Y) =    e^{{1\over
4h}|Y|^2}  \left ( J^{\Phi} _{Eh} W_E e^{{i\over \sqrt {h}} \Phi_S (
iY )} f \right ) (X)$$
for a.e. $X = (x ,
\xi)$ in $B(E) \times B(E)$.
Consequently:
$$e^{-{1\over 4h} |Y|^2} \  \rho_Y (  J_{ Eh}^{\Phi}  W _{E}
f   )    = <J^{\Phi} _{Eh} W_E e^{{i\over \sqrt {h}} \Phi_S ( - \eta
, y)} f , 1>  \ . $$
We deduce (4.30) since $1=J^{\Phi} _{Eh} W_E \Omega_K(E)$ and since $J^{\Phi} _{Eh}
W_E$ is an isometry.

\bigskip

{\it End of the proof of Theorem 4.10.}  From Proposition 4.11, we
have:
$$\rho _X \Big (\theta_{Eh} f \Big )  = e^{{1\over 4h} |X|^2}  <\Big ( J_{Eh}^{K} \Big
)^{-1} f , e^{- {i\over \sqrt {h}} \Phi_S (iX) }\Omega_K(E) >
_{{\cal H}(E)}$$
$$ = e^{{1\over 4h} |X|^2} < f , J_{Eh}^{K}e^{- {i\over \sqrt {h}} \Phi_S (iX) }
\Omega_K(E) > _{L^2( B(E) , \mu _{E,h}^K)}  \ . $$
We then apply Theorem 3.4 with the element $\Omega_K(E)$. Thus, we
derive (4.29) by direct computations.

\bigskip

Let us now prove that any function $F$ in $SB(E , h)$ is uniquely determined by its "restriction" to $Z_{\R}
\times Z_{\R}$.

\bigskip

{\bf Theorem 4.12.} {\it Suppose that $F$ in $SB(E , h)$ satisfies
$\rho_X (F) = 0$ for all $X$ in $Z_{\R} \times Z_{\R}$. Then $F
= 0$.

}

\bigskip

{\it  Proof.} Since $F$ is in $SB(E , h)$, there is
$f$ in ${\cal H}(E)$ such that $F = J_{Eh}^{\Phi} W_E f$. From Proposition 4.11, we have for each $X = (x , \xi)$ in $Z_{\R}
\times Z_{\R} = Z_{\bf C}$:
$$\rho_X (F) = e^{{1\over 4h} |X|^2} < f ,
e^{- {i\over \sqrt {h}} \Phi_S (iX) }\Omega_K(E) > _{{\cal H}(E)}  \ . $$

Under our hypothesis, the element $f$ is then orthogonal to all coherent states of ${\cal H}(E)$, namely, the elements
 $e^{- {i\over \sqrt {h}} \Phi_S (iX) }\Omega_K(E)$ ($X$
in $Z_{\bf C}$). According to proposition 3.1, the set of all these elements
is complete in  ${\cal H}(E)$. Therefore
$f=0$ and then $F =0$.

\bigskip

{\bf 5. Hybrid Weyl-anti-Wick  quantization.}

\bigskip

Let  $E$ be a finite subset of $\Gamma$ and  $G$  be a bounded
continuous function on $\R^E \times \R^E$, such that $\partial
_x^{\alpha} \partial _{\xi} ^{\beta}G$ is well-defined, bounded and
continuous on $\R^E \times \R^E$, for every multi-index $(\alpha ,
\beta)$ in $I_2(E)$. These hypotheses imply that the Weyl operator
$Op_h^{weyl}(G)$ is bounded in $L^2(\R^E)$. This is a result proved
by  H.O. Cordes for a formula similar to (1.1). It is also valid for
the Weyl formula (1.1) itself in view of the proof by A. Unterberger
[U2]. For all functions $\varphi$ and $\psi$ in ${\cal S}(\R^E)$, we
have:
$$<Op_h^{weyl, E}(G) \varphi , \psi> = (2  \pi h)^{-|E|} \int
_{(\R^E)^3} e^{{i\over h} (u-v)\cdot t} G \left ({u+v \over 2} , t
\right ) \varphi (v) \overline {\psi (u)} d\lambda _E(u) d\lambda
_E(v) d\lambda _E(t )  \ . \leqno (5.1)$$
In order to be self-contained, let us mention that (5.1) makes sense even without proving an $L^2$ estimate for this operator. To see this, it suffices to integrate by parts and to redefine the integral by
$$<Op_h^{weyl, E}(G) \varphi , \psi> = ...$$
$$ ... = (2  \pi h)^{-|E|} \int _{(\R^E)^3} e^{{i\over h} (u-v)\cdot \xi} G
\left ({u+v \over 2} , t \right )  \left ( 1 + {4|t|^2 \over h}
\right )^N (1 - (\partial _u - \partial _v)^2 )^N \Big ( \varphi (v)
\overline {\psi (u)}\Big )  d\mu $$
$$ d\mu = d\lambda _E(u) d\lambda _E(v) d\lambda _E(t)  \ . $$
This integral is convergent when choosing $N$ large enough, which gives a meaning to  (5.1).

\bigskip

To a subset $E$ and a function $G$ we may also associate an anti-Wick operator. In that case, the set $E$ may be infinite and the function $G$ is measurable and bounded on
$B(E) \times B(E)$. The operator $Op_h ^{AW, E} (G)$ is the bounded operator in ${\cal H}(E)$ satisfying,
$$ < Oph^{AW , E}(G) f , g > _{{\cal H}(E)} = \int _{B(E) \times
B(E)} G(X)  (J^\Phi_{Eh} W_E f ) (X) \ \overline { (J^\Phi_{Eh} W_E g)
(X)} d\mu_{E,h}^{\Phi}(X)\leqno (5.2) $$
for every
$f$ and $g$ in ${\cal H}(E)$, where the Segal isomorphism $
J^\Phi_{Eh} $ has been defined in Section 3.B.

\bigskip

We shall now define a hybrid operator $Op_h^{hyb, E}(F)$ bounded in
${\cal H}(\Gamma)$, for every finite subset $E$ of $\Gamma$ and
every function
 $F$ on $B(\Gamma) \times
B(\Gamma)$. It suffices to define the scalar product
$< OP_h^{hyb , E} (F) f , g>$ for all $f$ and $g$ in ${\cal
H}(\Gamma)$. For the sake of clarity, let us first begin with the case when:
$$ f = f_E \otimes f_{E^c} \hskip 2cm g = g_E \otimes g_{E^c}\ , \leqno (5.3)$$
with $f_E$ and $g_E$ in ${\cal H} (E)$,  $f_{E^c}$ and $g_{E^c}$ in
${\cal H} (E^c)$. The running variable in $B(\Gamma) \times
B(\Gamma)$ may be written as $( X_E , X_{E^c})$. For every $X_{E^c}$
in $B(E^c) \times B(E^c)$,   $F_{X_{E^c}}$ stands for the function
defined on $\R^E \times \R^E$ by:
$$F_{X_{E^c}} (X_E) = F( X_E , X_{E^c}) \hskip 2cm X_E \in \R^E \times
\R^E \ . \leqno (5.4)$$
When  $f$ and $g$ are written as in (5.3), the hybrid operator shall be defined as,
at least formally:
$$< OP_h^{hyb , E} (F) f , g> _{{\cal H}(\Gamma)} = < Op_h^{AW ,
E^c} (\Phi) f_{E^c} , g_{E^c} > _{{\cal H}(E^c)}\ , \leqno (5.5)$$
where $\Phi$ is the function defined on $B(E^c) \times B(E^c)$
by:
$$ \Phi (X_{E^c}) = < Op_h^{weyl , E} (F_{X_{E^c}}) f_E , g_E >
 _{{\cal H}(E)}\ , \leqno (5.6)$$
with $F_{X_{E^c}}$ given by (5.4).

\bigskip

Next, we consider $f$ and $g$ in ${\cal H}(\Gamma)$, not necessarily
written as in (5.3). We shall use the coherent states  $\Psi_{X_E,
h}$ defined in (4.1) with finite $E$. For every $X_E$ in $\R^E
\times \R^E$ ($E$ finite),
 $i^{\star}_{X_E}$ denotes the mapping from ${\cal H}(\Gamma)$ into
${\cal H}(E^c)$  satisfying:
$$ < i^{\star}_{X_E} f , g > _{{\cal H}(E^c)} = <f ,\phi_{X_{E h}}
\otimes g >_{{\cal H}(\Gamma)}\ , \leqno (5.7)$$
with $\phi_{X_{E h}}$ defined in (4.19).
Following (4.2), we may write, at least formally
$$ f = (2 \pi h)^{-|E|} \int _{\R^E \times \R^E} \phi_{X_{E  h}}
\otimes i^{\star}_{X_E} f \ d\lambda _E(X_E) \leqno (5.8)$$
for any
$f$ in ${\cal H}(\Gamma)$, with a finite $E$ and where $\lambda _E$ is the Lebesgue measure.

\bigskip

The expansion (5.8) may be used to verify that the definition of the
hybrid operator in (5.5),(5.6) is a particular case of  Definition
5.1 below. In view of Definitions (5.1) and (5.2) of the
 Weyl and anti-Wick operators, we are led  in the general case to the
 following definition.

\bigskip

{\bf Definition 5.1.} {\it For every bounded continuous function $F$ on $B (\Gamma) \times B(\Gamma) $ and for any finite subset
$E$ of $\Gamma$, we denote by $Op_h^{hyb, E}(F)$
the operator in ${\cal H}(\Gamma)$ such that, for all $f$ and $g$
in ${\cal H}(\Gamma)$:
$$ < Op_h^{hyb, E}(F) ( f ) , g  > = ... \leqno (5.9) $$
$$ =(2\pi h)^{-2|E|}   \int _{\Delta (E) }
<Op_h^{weyl, E}(F_{Z_{E^c}}) \Psi^E _{X_E h} , \Psi^E_{Y_E h} >   \
(J_{E^c h}^{\Phi}W_{E^c } i _{X_E h}^{\star} f )(Z_{E^c})  \
\overline { (J_{E^c h}^{\Phi}W_{E^c } i _{Y_E h}^{\star} g )
(Z_{E^c})} ... $$
$$ ... d\lambda _E(X_E) d\lambda _E ( Y_E)  d\mu _{E^c,h}^{\Phi} (Z_{E^c})\ ,$$
where $\Delta (E)=(\R^2)^E \times (\R^2)^E \times (B (E^c)\times B
(E^c))$.  The scalar product in the above right hand-side (resp. left hand-side) is the scalar product in ${{\cal H}(E)}$ (resp.  ${\cal H}(\Gamma)$).
Te Segal isomorphism
$J_{E^c h}^{\Phi}$ has been diefined in Section 3.B.
 When  $E$ is empty, $Op_h^{hyb, E}$ is defined as the anti-Wick operator.
}

\bigskip

We now state the main Theorems.

\bigskip

{\bf Theorem 5.2.} {\it  Set $0< h \leq 1$. Let $F$ be a function
satisfying the hypothesis $H_2(M, \varepsilon)$ in Definition 1.3
with
   a real number $M>0$ and with a family $\varepsilon = (\varepsilon _j)_{(j\in \Gamma)}$  of nonnegative real numbers. Let  $\Lambda$ be a finite subset of  $\Gamma$.
   Then, the operator $Op_h^{(hyb , \Lambda)}(F)$ in Definition 5.1 is bounded in ${\cal H}(\Gamma)$ and we have:
$$ \Vert  Op_h^{hyb, \Lambda} ( F) \Vert _{{\cal L}({\cal H})} \leq M \prod
_{j\in \Lambda} ( 1 + 225 \pi K_2 \sqrt {h} \varepsilon _j)\ , \leqno
(5.10)$$
where $K_2 = \sup _{j\in \Gamma}\ \max (1 , \varepsilon_j^3)$.  Under hypothesis $H_4(M, \varepsilon)$, we have:
$$ \Vert  Op_h^{hyb, \Lambda} ( F) \Vert _{{\cal L}({\cal H})} \leq M \prod
_{j\in \Lambda} ( 1 + 225 \pi K_4 h \varepsilon _j^2) \ ,\leqno
(5.11)$$
where $K_4= \sup _{j\in \Gamma}\ \max (1 , \varepsilon_j^6)$.

}

\bigskip

{\bf Theorem 5.3.} {\it Let $0 <h \leq 1$. Suppose that $F$ is a
function as in Theorem 5.2 and assume that  $\Lambda$ and $\Lambda'$
are two finite subsets of $\Gamma$ such that $\Lambda \subset
\Lambda'$. Then, we have:
$$ \Vert  Op_h^{hyb , \Lambda'} (F)- Op_h^{hyb , \Lambda} (F)
 \Vert _{{\cal L}({\cal H})} \leq M 225 \pi
K_2 \sqrt {h}\left [  \sum _{j\in \Lambda' \setminus \Lambda }
\varepsilon _j \right ] \ \prod _{k\in \Lambda'} (1 + 225 \pi K_2
\sqrt {h} \varepsilon_k)\ . \leqno (5.12)$$
If $H_4(M , \varepsilon)$ is satisfied, then
$$ \Vert  Op_h^{hyb , \Lambda'} (F)- Op_h^{hyb , \Lambda} (F)
 \Vert _{{\cal L}({\cal H})} \leq M 225 \pi
K_4 h  \left [  \sum _{j\in \Lambda' \setminus \Lambda } \varepsilon
_j^2 \right ] \ \prod _{k\in \Lambda '} (1 + 225 \pi K_4 h
\varepsilon_k^2)\ . \leqno (5.13)$$

 }

\bigskip

{\bf Theorem 5.4.} {\it  Suppose that $F$ is a function as in
Theorem 5.2.  and assume that $(\Lambda_n)$ is an increasing
sequence of finite subsets of $\Gamma$ with union $\Gamma$. Then,
the sequence of operators $( Op_h^{hyb , \Lambda_n} (F))_{(n\in
\N)}$ is a Cauchy sequence in  ${\cal L}({\cal H}(\Gamma))$. Its
limit is denoted by $Op_h^{weyl} (F)$. It is independent of the
chosen sequence and satisfies, if  $0 <h \leq 1$:
$$\Vert  Op_h^{weyl} (F) \Vert _{{\cal L}({\cal H})}
\leq M \prod _{j\in \Gamma } (1 + 225 \pi K_2 \sqrt {h} \varepsilon
_j)\ . \leqno (5.14)$$
When the hypothesis $H_4(M , \varepsilon)$ is verified, we have:
$$\Vert  Op_h^{weyl} (F) \Vert _{{\cal L}({\cal H})}
\leq M \prod _{j\in \Gamma } (1 + 225 \pi K_4 h \varepsilon _j^2
)  \ . \leqno (5.15)$$

 }

Theorem 5.4 is a consequence of Theorem 5.2 and Theorem 5.3. If
$(\varepsilon_j)_{(j\in \Gamma )}$ is a summable family then
$$ \lim _{n\rightarrow \infty }\sum _{j\notin \Lambda_n} \varepsilon _j
= 0\ , $$
where $(\Lambda_n)$ is a sequence as in Theorem 5.4 and the infinite
product in (5.14) is convergent. The same holds true if
$(\varepsilon_j^2)_{(j\in \Gamma)}$ is summable.

\bigskip

{\bf 6. Changing the subset $E$.}

\bigskip

We have to study the relation between $Op_h^{hyb, E_1}(F)$
and $Op_h^{hyb, E_2}(F)$
when the two subsets  $E_1$ and $E_2$  of $\Gamma$ are finite and verify
$E_1 \subset E_2$ and when the function $F$ satisfies our hypotheses.

\bigskip

\bigskip

{\bf Proposition 6.1.} {\it For every finite subsets  $E_1$ and
$E_2$  of $\Gamma$ such that $E_1 \subset E_2$ and for all functions
$F$ satisfying the hypothesis in Theorem 5.2, we have:
$$Op_h^{(hyb , E_1)}( F) = Op_h^{(hyb , E_2)}
\left (e^{{h\over 4} \Delta _{E_2 \setminus E_1}} F \right )\ , \leqno
(6.1)$$
where:
$$ \Delta _{E_2 \setminus E_1} = \sum _{j\in E_2 \setminus E_1}
{\partial ^2 \over \partial x_j^2} + {\partial ^2 \over \partial
\xi_j^2}\ .\leqno (6.2)$$

}

\bigskip

We make use of the following notations throughout the proof. We set
$D = E_2 \setminus E_1$ and $S = (E_2)^c$. Thus, the variable in
$\R^{E_2}$ may be written as $(x_{E_1}, x_D)$ and the one in $B(
E_1^c)$ is denoted by  $(x_D , x_S)$. The proof is obtained by
combining  Definition 5.1 with   $E_1$ and $E_2$, together with the
following two  lemmas.

\bigskip

{\bf Lemma 6.2.}   {\it For every continuous function $F$ on
$B(\Gamma) \times B(\Gamma)$ satisfying the hypotheses in Theorem
5.2, for any $X_{E_2} = ( X_{E_1} , X_D)$ and $Y_{E_2} = ( Y_{E_1} ,
Y_D)$ in $(\R^{E_2})^2$,  we have:
$$ <Op_h^{(weyl, E_2)}
\Big (e^{{h\over 4} \Delta _{D}} F_{Z_S}\Big )  \Psi_{X_{E_1}, X_D,
h} , \Psi_{Y_{E_1}, Y_D,  h}  > = ... \leqno (6.3)$$
$$ ... =  (2 \pi h)^{-|D|} \int _{(\R^D )^2}   <Op_h^{(weyl, E_1)}(F_{Z_D, Z_S}) \Psi_{X_{E_1} h} , \Psi_{Y_{E_1}
h}  > < \Psi_{X_{D} h} , \Psi_{Z_{D} h}> < \Psi_{Z_{D} h} ,
\Psi_{Y_{D} h}> d\lambda_D (Z_D) \ .$$
The scalar product in the above left hand-side  is the scalar product of $L^2(\R^{E_2})$ and in the right hand-side, it is the scalar product of  $L^2(\R^{E_1})$ (for the first one) and  $L^2(\R^{D})$
(for the last two ones).
Also,  $F_{Z_D, Z_S}$ stands for the function defined on  $ (\R^{E_1})^2$ by $F_{Z_D, Z_S}(Z_{E_1}) = F(Z_{E_1}
, Z_D, Z_S)$, for all $Z_{E_1}$ in $(\R^{E_1})^2$.

}
\bigskip

{\it Proof.} The left hand-side of (6.3) is expressed using the
Weyl calculus definition (5.1), replacing $E$ by
$E_2 = E_1 \cup D$, $G$ being replaced by $e^{{h\over 4} \Delta _{D}} F_{Z_S} $,
$u$, $v$ and $t$ being replaced by $(u_{E_1} , u_D)$, $(v_{E_1} , v_D)$
and $(t_{E_1} , t_D)$, $\varphi (v)$ by $\Psi _{X_{E_1} , h}
(v_{E_1}) \Psi _{X_{D} , h} (v_{D})$, and $\psi (u)$ by $\Psi
_{Y_{E_1} , h} (u_{E_1}) \Psi _{Y_{D} , h} (u_{D})$. Next, we write:
$$(e^{{h\over 4} \Delta _{D}} F_{Z_S} ) \left ( { u_{E_1} +
v_{E_1} \over 2} , {u_D + v_D\over 2}, t_{E_1} , t_D \right ) =
...$$
$$ ... =
(\pi h) ^{-|D|} \int _{\R^D \times \R^D} e^{ - {1 \over h} \left
|{u_D + v_D\over 2} - z_D \right |^2-{1 \over h} | t_D - \zeta _D
|^2} F_{Z_S}  \left ( { u_{E_1} + v_{E_1} \over 2}  , z_D ,t_{E_1} ,
\zeta _D \right ) dz_D d\zeta_D\ .$$
Then, setting $Z_D = (z_D , \zeta_D)$ we remark that:
$$ (\pi h)^{-|D|} \int _{\R^D} e^{ - {1 \over h} \left
|{u_D + v_D\over 2} - z_D \right |^2-{1 \over h} | t_D - \zeta _D
|^2} e^{{i\over h} (u_D - v_D) t_D} dt_D = \Psi_{Z_D , h} (u_D)
\overline {\Psi_{Z_D , h} (v_D) }\ . $$
According to Definition 5.1 with $E$ now replaced by $E_1$ and $G$
replaced by $F_{Z_D , Z_S}$, it leads to the right hand-side of (6.3).

\bigskip

{\bf Lemma 6.3.} {\it For any $f$ in ${\cal H}(\Gamma)$, for any
$X_{E_1}$ in $(\R^{E_1})^2$, for any $Z_D$ in $(\R^D)^2$, for a.e.
$Z_S$ in $B(S) \times B(S)$, we have:
$$e^{-{|Z_D|^2 \over 4h}}  \Big ( J_{D \cup S , h}^{\Phi} W_{D\cup S} i_{X{E_1}}^{\star} f
\Big ) (Z_D , Z_S) = ... \leqno (6.4)$$
$$ ... = (2 \pi h)^{-|D|}  \int _{(\R^D)^2}  < \Psi_{X_D , h} , \Psi_{Z _D , h} >
_{L^2(\R^d)} \ \Big ( J_{ S , h}^{\Phi} W_{ S} i_{X_{E_1},
X_D}^{\star} f \Big ) ( Z_S) d\lambda (X_D)\ .
$$

}

\bigskip

{\it Proof.} It is sufficient to prove (6.4) when $f =
f_{E_1}\otimes f_{D} \otimes f_S$, with $f_{E_1}$ in ${\cal
H}(E_1)$, $f_D$ in ${\cal H}(D)$, $f_S$ in ${\cal H}(S)$. Then, from (5.7):
$$  i _{X_{E_1}} f = < f_{E_1} , \varphi _{X_{E_1}, h} >_{{\cal H}(E_1)}  f_D \otimes
f_S\ . $$
Consequently, from (3.18) and  (4.7):
$$\Big ( J_{D \cup S , h}^{\Phi} W_{D\cup S} i_{X{E_1}}^{\star} f
\Big ) (Z_D , Z_S) = < f_{E_1} , \varphi _{X_{E_1}, h} >_{{\cal
H}(E_1)}  \Big (J^{\Phi} _{Dh} W_D f_D \Big ) (Z_D) \Big (J^{\Phi}
_{S h} W_S f_S \Big ) (Z_S)\ .$$
Since $D$ is finite, we have, using Proposition 4.11:
$$e^{-{|Z_D|^2 \over 4h}}  \Big (J^{\Phi} _{Dh} W_D f_D \Big ) (Z_D) =
< f_D , \varphi _{Z_D , h} >_{{\cal H} (D)}  $$
and in view of (4.2):
$$ ... = (2 \pi h)^{-|D|}  \int _{\R^D \times \R^D}
< f_D , \varphi _{X_D , h} >_{{\cal H} (D)}  < \Psi_{X_D , h} , \Psi
_{Z_D , h} > _{L^2(\R^d)} d\lambda_D ( X_D)\ .
$$
According to (5.7), we also have:
$$< f_{E_1} , \varphi _{X_{E_1}, h} >_{{\cal H}(E_1)}
< f_{D} , \varphi _{X_{D}, h} >_{{\cal H}(D)} f_S =
i^{\star}_{X_{E_1} , X_D} f\ .$$
Equality (6.4) then follows when $f = f_{E_1}\otimes f_{D}
\otimes f_S$. In view of (3.8), it also holds true in the general case using linearity and continuity of the two hand-sides.

\bigskip

{\it End of the proof of Proposition 6.1.}  With the above notations, it suffices to show:
$$ <Op_h^{(hyb , E_2)}
(e^{{h\over 4} \Delta _{D}} F) f , g> = <Op_h^{(hyb , E_1)}( F) f ,
g >   \leqno (6.5)
$$for all $f$ and $g$
in ${\cal H}(\Gamma)$.
The left hand-side is rewritten using Definition 5.1 with the set
$E_2 = E_1 \cup D$ and the function $e^{{h\over 4} \Delta _{D}} F$.
The $L^2(\R^{E_2})$ scalar product in the new expression is
rewritten using Lemma 6.2. Then,  Lemma 6.3 is applied twice with
$f$, $X_{E_1}$ and $Z_D$ and once with
 $g$, $Y_{E_1}$ and $Z_D$. Next we observe, from (2.8):
$$ d\mu ^{\Phi} _{ D \cup S, h} (Z_D , Z_S) = (2 \pi h)^{-|D| }
e^{ -{|Z_D|^2 \over 2h}} d\lambda _D(Z_D) d\mu ^{\Phi} _{  S, h}
(Z_S)\ . $$
In this new expression, it appears the right hand-side of (6.5)
using Definition 5.1 with the set $E_1$ and with the function $F$.
This proves (6.5) and the proof of the proposition is completed.
\bigskip

{\bf 7. First reductions.}

\bigskip

Let $\Lambda$ be a finite subset of $\Gamma$. Denoting by
$I_{\Lambda}$ the identity operator in the space of bounded continuous functions on $\R^{\Lambda}\times \R^{\Lambda}$, we have:
$$ I_{\Lambda} = \sum _{E \subseteq \Lambda}
 \prod _{j\in E} ( I - e^{{h\over 4}\Delta_j} )
\prod _{j\in \Lambda \setminus  E} e^{{h\over 4} \Delta _j}\ . \leqno
(7.1)$$
In the above equality, the sum is running
over all subsets of $\Lambda$ including the empty subset and
$$\Delta _j =  {\partial^2 \over \partial x_j^2} +
{\partial^2 \over \partial \xi_j^2}\ .  \leqno (7.2)$$

\bigskip

The operator $T_h (E)$ stands for:
$$T_h (E) = \prod _{j\in E} ( I - e^{{h\over 4}\Delta_j} ) \ \leqno (7.3)$$
for every finite subset $E$ of $\Gamma$ and every $h>0$. If  $E$ is
empty then  $T_h (E) = I$.  From equality (7.1) and  Proposition 6.1
we obtain the following identity.

\bigskip

{\bf Proposition 7.1.} { \it For every function $F$ satisfying
hypothesis $H_2(M , \varepsilon)$ in Definition 1.3 and for any
finite subset $\Lambda$ in $\Gamma$, we have:
$$Op_h^{(hyb , \Lambda)} (F) = \sum _{E \subseteq \Lambda} Op_h^{hyb, E}
(T_h(E) F)\ . \leqno (7.4)$$
Again, the sum is running over all subsets of $\Lambda$ including the empty subset. We have
$$Op_h^{hyb, \emptyset} (F)= Op_h^{AW, \Gamma } (F)\ .$$
For every finite subset $\Lambda$ and $\Lambda'$ with $\Lambda
\subset \Lambda '$, we have:
$$Op_h^{hyb, \Lambda ' } (F) - Op_h^{hyb , \Lambda } (F)=
\sum _{E \in P( \Lambda, \Lambda ')} Op_h^{(hyb , E)} (T_h(E) F)\ ,
\leqno (7.5)$$
where $P( \Lambda, \Lambda ')$ is the set of all the  $E \subseteq
\Lambda ' $, which are not included in $\Lambda$ and, 
in particular,  not empty.

}

\bigskip

The most technical part of this work is the following proposition.
It will be proved in Section 8.

\bigskip

{\bf Proposition 7.2.} {\it For every function $F$ verifying
hypothesis $H_2(M , \varepsilon )$ in Definition 1.3 and for any
finite subset  $E$ of $\Gamma$, we have
$$ \Vert  Op_h^{(hyb , E)} (T_h(E) F) \Vert _{{\cal L}({\cal H}(\Gamma))}
\leq M (225 \pi K_2 \sqrt {h})^{|E|} \prod _{j\in E}  \varepsilon
_j\ , $$
where $K_2 = \sup _{j\in \Gamma} \ \max (1 , \varepsilon_j^3)$. If
$E$ is empty then the norm is bounded by $M$. If hypothesis
$H_4(M , \varepsilon )$ is satisfied then:
$$ \Vert  Op_h^{(hyb , E)} (T_h(E) F) \Vert _{{\cal L}({\cal H}(\Gamma))}
\leq M (225 \pi K_4 h)^{|E|} \prod _{j\in E}  \varepsilon _j^2\ , $$
where $K_4 = \sup _{j\in \Gamma} \ \max (1 , \varepsilon_j^6)$.

}

\bigskip

Theorem 5.2 and Theorem 5.3 are easily deduced from the two above
propositions.

\bigskip

{\it Proof of Theorem 5.2.} If $H_2(M , \varepsilon)$ is verified,
it is seen from Proposition 7.1 (point (7.4)) and Proposition 7.2
that:
$$ \Vert  Op_h^{(\Lambda)} (F) \Vert _{{\cal L}({\cal H})} \leq M
\sum _{E \subseteq \Lambda} (225 \pi K_2  \sqrt {h})^{|E|}
\prod_{j\in E} \varepsilon_j = M \prod _{j \in \Lambda}  (1 + 225
\pi K_2 \sqrt {h} \varepsilon _j)\ ,$$
where $K_2 = \sup _{j\in \Gamma}\ \max (1 , \varepsilon_j^3)$. When
$H_4(M , \varepsilon)$ is verified, we have:
$$ \Vert  Op_h^{(\Lambda)} (F) \Vert _{{\cal L}({\cal H})} \leq M
\sum _{E \subseteq \Lambda} (225 \pi K_4 h)^{|E|} \prod_{j\in E}
\varepsilon_j^2 = M \prod _{j \in \Lambda}  (1 + 225 \pi K_4  h
\varepsilon _j^2)\ , $$
where $K_4= \sup _{j\in \Gamma}\ \max (1 , \varepsilon_j^6)$.

\bigskip

{\it Proof of Theorem 5.3.} According to Proposition 7.1 (point
(7.5)) and Proposition 7.2,  when $\Lambda $ and $\Lambda '$ are two
finite subsets of $\Gamma $ with $\Lambda \subset \Lambda '$ and if
$H_2(M , \varepsilon)$ is satisfied, we have:
$$ \Vert  Op_h^{(\Lambda')} (F)- Op_h^{(\Lambda)} (F) \Vert _{{\cal L}({\cal H})} \leq M
\sum _{E \in P( \Lambda, \Lambda ') } (225 \pi K_2 \sqrt {h})^{|E|}
\prod_{j\in E} \varepsilon_j\ ,$$
where $P( \Lambda, \Lambda ')$ is the set of  $E$ in
$\Lambda '$ not being included in $\Lambda$.
Inequality (5.12) then follows.

\bigskip

The norm $N(f)$ defined below for every $f$ in ${\cal H}(\Gamma)$ and for any finite subset
$E$ of $\Gamma$ shall be involved in the proof of Proposition 7.2 and then also in our main results. Set
$$ N_E(f)^2 =  (2 \pi h)^{-|E|} \int _{(\R^E)^2 \times B (E^c)^2}
|(J_{E^c h}W_{E^c} i _{X_E h}^{\star} f )(Z_{E^c}) |^2 dX_E d\mu
_{E^c, h} ^{\Phi} (Z_{E^c})\ . \leqno (7.8) $$
The next proposition will be useful.

\bigskip

{\bf Proposition 7.3.} {\it For every $f$ in ${\cal H}(\Gamma)$ and for any finite subset $E$ of $\Gamma$, we have $N_E(f) = \Vert f
\Vert _{{\cal H}(\Gamma)}$.

}

\bigskip

{\it Proof.} Since $J_{E^c h}$ is an isometric isomorphism between ${\cal H}_{\Phi}(E^c)$ and  $L^2( B(E^c)^2 , \mu
_{E^c, h} ^{\Phi})$ and since $W_{E^c}$ is a partial isometry from ${\cal H}(E^c)$ to ${\cal H}_{\Phi}(E^c)$,  then we have:
$$ \int _{ B (E^c)\times B(E^c)} |(J_{E^c h}W_{E^c } i _{X_E h}^{\star} f )(Z_{E^c})
|^2 d\mu _{E^c,h} (Z_{E^c})  =  \Vert W_{E^c } i _{X_E h}^{\star} f
\Vert _{{\cal H}_{\Phi}(E)}^2 =  \Vert i _{X_E h}^{\star} f \Vert
_{{\cal H} (E^c) }^2 \leqno (7.9)$$
for all $X_E$ in $(\R^2)^E$. We shall now show that:
$$ (2 \pi h)^{-|E|} \int _{(\R^E)^2} \Vert i _{X_E h}^{\star} f \Vert
_{{\cal H} (E^c) }^2 = \Vert f \Vert _{{\cal H}}^2 \ .$$
More generally, we shall prove that:
$$ < f , g> _{{\cal H}}= (2 \pi h)^{-|E|}
\int _{(\R^E)^2} < i_{X_E}^{\star} f , i_{X_E}^{\star} g
> _{{\cal H}(E^c)}dX_E \leqno (7.10) $$
 for all $f$ and $g$
in  ${\cal H}={\cal H}(\Gamma) $.
We first prove  (7.10) when $f =
f_E \otimes f_{E^c}$ and $g = g_E \otimes g_{E^c}$ with $f_E$ and
$g_E$ in ${\cal H} (E)$, $f_{E^c}$ and $g_{E^c}$ in ${\cal H}
(E^c)$. In that case, we have:
$$(2 \pi h)^{-|E|} \int _{(\R^E)^2} <
i_{X_E}^{\star} f , i_{X_E}^{\star} g > _{{\cal H}(E^c)}dX_E = ...$$
$$ ... = (2 \pi h)^{-|E|} \int _{(\R^E)^2} < \widetilde J_{Eh}^Kf_E , \Psi_{X_E h}> < \Psi_{X_E
h} , \widetilde J_{Eh}^K g_E> <f_{E^c} , g_{E^c}> _{{\cal H}(E^c)}
dX_E$$
$$ =<f_{E} , g_{E}>_{{\cal H}(E)}   <f_{E^c} , g_{E^c}>_{{\cal H}(E^c)}
= < f , g>_{\cal H}$$
We have here used (4.2). Equality (7.10) then follows in the general case since both terms are continuous bilinear mappings on  ${\cal H}(\Gamma)$ and since ${\cal
H}(\Gamma) $ is the completion of the tensor product ${\cal H}(E)
\otimes {\cal H}(E^c)$.

\bigskip

{\bf 8. Proof of Proposition 7.2.}

\bigskip

We shall first give a bound on the scalar product appearing in the integral (5.9) defining the hybrid quantization, when $F$ is replaced
by $T_h(E) F$ with $F$ satisfying hypothesis
$H_2(M, \varepsilon)$ and where
$T_h(E)$ is defined in (7.3).

\bigskip

{\bf Proposition 8.1.} {\it If $F$ verifies hypothesis $H_2(M,
\varepsilon)$ in Definition 1.3  with a constant $M>0$ and a
summable family $(\varepsilon_j)_{(j\in \Gamma)}$, if $E$ is a
finite subset of $\Gamma$, if $0 < h \leq1$ and if $Z_{E^c}$ is in
$B( E^c) \times B(E^c)$, then:
$$
| <Op_h^{weyl, E}(T_h(E) F_{Z_{E^c}}) \Psi^E _{X_E h} , \Psi^E_{Y_E
h} >  | \leq ... \leqno (8.1) $$
$$ ... \leq M (450 K_2 \sqrt {h})^{|E|}\prod _{j\in E} \varepsilon_j \left ( 1 +
{|x_j - y_j|^2 \over h} \right )^{-1} \left ( 1 + {|\xi_j -
\eta_j|^2 \over h} \right )^{-1} \ , $$
where $K_2 = \sup _{j\in \Gamma} \max ( 1 , \varepsilon_j^3)$. When the hypothesis $H_4(M, \varepsilon)$  is satisfied, we have:
$$
| <Op_h^{weyl, E}(T_h(E) F_{Z_{E^c}}) \Psi^E _{X_E h} , \Psi^E_{Y_E
h} >  | \leq ...$$
$$ ... \leq M (450 K_4 h)^{|E|}\prod _{j\in E} \varepsilon_j^2 \left ( 1 +
{|x_j - y_j|^2 \over h} \right )^{-1} \left ( 1 + {|\xi_j -
\eta_j|^2 \over h} \right )^{-1} \ , $$
where $K_4 = \sup _{j\in \Gamma} \max ( 1 , \varepsilon_j^6)$.

}

\bigskip

Proposition 8.1 will rely on propositions 8.2, 8.3 and 8.4 below. The first one is concerned with an integral expression of the left hand-side of (8.1).

\bigskip

{\bf Proposition 8.2.} {\it For any finite subset  $E$
of $\Gamma$, for every bounded continuous function  $G$ on $\R^E
\times \R^E$, for all $X_E$ and $Y_E$ in $(\R^E)^2 $, we have:
$$ <Op_h^{weyl,  E}(G) \Psi^E _{X_E h} , \Psi^E_{Y_E h} >  = (\pi h)^{-|E|}
 \int  _{(\R^E)^2} G (Z_E ) e^{ - {1 \over h}
\left | Z_E - {X_E + Y_E \over 2}\right |^2 } \ e^{{i\over h}
\varphi (X_E , Y_E , Z_E) } d\lambda _E (Z_E) \leqno (8.2) $$
setting $X_E = (x_E , \xi_E)$, $Y_E = (y_E , \eta_E)$,
$Z_E = (z_E , \zeta_E)$ and
$$\varphi (X_E , Y_E , Z_E) = z_E \cdot (\xi_E - \eta _E) -
 \zeta_E \cdot ( x_E - y_E) + {1\over 2} ( x_E
\cdot \eta_E - y_E  \cdot \xi_E) \ .\leqno (8.3)$$

 }

\bigskip

This proposition is derived by Unterberger [U2] (formula (1.3)). In
[U2], the integral kernel in (8.2) is defined as the Wigner function
associated with the two coherent states $\Psi^E _{X_E h}$ and
$\Psi^E _{Y_E h}$, taken at the point $Z_E$. It is expressed in
(1.4) of [U2] and direct computations give (8.2) and (8.3).

\bigskip

{\bf Proposition 8.3.} {\it For every function $G$ satisfying
hypothesis $H_2(M , \varepsilon)$ in Definition 1.3, for each
$Z_{E^c}$ in $B(E^c) \times B(E^c)$ and for each $h>0$:
$$ \prod _{j\in E} \left ( 1 + {|x_j - y_j|^2 \over h} \right )
\left ( 1 + {|\xi_j - \eta_j|^2 \over h} \right ) | <Op_h^{weyl,
E}(G_{Z_{E^c}} ) \Psi^E _{X_E h} , \Psi^E_{Y_E h} > | \leq ...
\leqno (8.4)$$
$$...  \leq C^{|E|} \sum _{ (\alpha , \beta)  \in
I_2(E) }h^{(|\alpha|+|\beta|)/2}
 \Vert
 \partial _z^{\alpha} \partial _{\zeta}^{\beta} G_{Z_{E^c}}  \Vert
 _{L^{\infty}(\R^E \times \R^E)}\ , $$
where $C = 25$ and
  $I_m(E) = \{ 0, 1 , ... , m \}^E \times  \{ 0, 1 , ... , m
 \}^E$  for all $m\geq 1$,.

}

 \bigskip

{\it Proof.}  Integrating by parts (8.2) with $G$ replaced by $G_{Z_{E^c}}$  yields:
$$ |x_j - y_j|^2 <Op_h^{weyl,  E}(G_{Z_{E^c}} ) \Psi^E _{X_E h} , \Psi^E_{Y_E h} > = ...$$
$$ =  - h^2 (\pi h)^{-|E|} \int  _{(\R^E)^2} e^{{i\over h} \varphi (X_E
, Y_E , Z_E) }  {\partial^2 \over \partial \zeta_j^2}  \left [ G
(Z_E , Z_{E^c} ) e^{ - {1 \over h} \left | Z_E - {X_E + Y_E \over
2}\right |^2 } \right ]d\lambda (Z_E)
$$
for   $j\in E$, where $\varphi$ is defined in (8.3). Iterating this process, we obtain:
$$ \prod _{j\in E} \left ( 1 + {|x_j - y_j|^2 \over h} \right )
\left ( 1 + {|\xi_j - \eta_j|^2 \over h} \right ) |<Op_h^{weyl,
E}(G_{Z_{E^c}} ) \Psi^E _{X_E h} , \Psi^E_{Y_E h} >|
 \leq  ... $$
$$ \leq (\pi h)^{-|E|} \int  _{(\R^E)^2} \left | H(X_E , Y_E , Z_E, Z_{E^c})
  \right |d\lambda (Z_E)\ ,
$$
where:
$$H(X_E , Y_E , Z_E, Z_{E^c} )  = \prod _{j\in E}
\left ( 1 - h  {\partial^2 \over \partial z_j^2} \right ) \left ( 1
- h  {\partial^2 \over \partial \zeta_j^2} \right ) \left [ G (Z_E ,
Z_{E^c} ) e^{ - {1 \over h} \left | Z_E - {X_E + Y_E \over 2}\right
|^2 } \right ]\ . $$
Clearly:
$$H(X_E , Y_E , Z_E , Z_{E^c})  =e^{ - {1 \over h} \left | Z_E - {X_E + Y_E \over
2}\right |^2 } \prod _{j\in E}L_{z_ j} L_{\zeta _j }G (Z_E , Z_{E^c}) \ , $$
where we use the notation:
$$L_{z _j} = \sum _{k=0}^3 p_k \left ( h^{-1/2} \left ( z_j - { x_j +
y_j \over 2} \right ) \right ) h^{k/2} {\partial^k \over \partial
z_j^k}\ , $$
where
$$ p_0(x) =  3 - 4 x^2 \hskip 2 cm p_1(x) = 4x \hskip 2cm p_2(x) =
-1$$
and with similar notations for $L_{\zeta _j}$. We may write:
$$\prod _{j\in E}L_{z _j} L_{\zeta _j} = \sum _{(\alpha , \beta ) \in
I_2(E)  } h^{(|\alpha| +|\beta|)/2} A_{\alpha \beta} \left (
h^{-1/2} \left ( Z_E - { X_E + Y_E \over 2} \right ) \right )
\partial _z^{\alpha}
\partial _{\zeta}^{\beta} $$
with:
$$ A_{\alpha \beta}(z , \zeta) = \prod _{j\in E} p_{\alpha _j }
(z_j) p_{\beta _j } (\zeta_j)\ .$$
Let $C>0$ be such that:
$$ C^{1\over 2} \geq  \pi ^{-1/2} \int _{\R} |p_k(x)| e^{-x^2} dx
\hskip 2cm k=0, 1 , 2\ . \leqno (8.5)$$
Inequality (8.4) then holds true. One may choose $C=25$ in order to satisfy (8.5).

 \bigskip

Now, we shall apply
 Proposition 8.3 setting $G =
T_h(E)F$ where $T_h(E)$ is defined in (7.3) and $F$ verifies hypothesis $H_2(M, \varepsilon)$. Replacing $G$ by
$T_h(E)F$, the constants in the right hand-side of (8.4) will be improved.

\bigskip

{\bf Proposition 8.4.} {\it Suppose $0 \leq h < 1$ and assume that $F$
satisfies hypothesis $H_2(M, \varepsilon)$. For any $Z_{E^c}$  $B(E^c) \times B(E^c)$, we have:
$$ \sup _{(\alpha , \beta)   \in
I_2(E) } h^{(|\alpha|+|\beta|)/2}
  \Vert
 \partial _z^{\alpha} \partial _{\zeta}^{\beta} T_h(E) F_{Z_{E^c}}  \Vert
 _{L^{\infty}(\R^{E}\times \R^E)} \leq 2^{|E|}  \sup _{(\alpha , \beta)\in \widetilde
 I_2(E)
 } h^{(|\alpha|+|\beta|)/2}
  \Vert
 \partial _z^{\alpha} \partial _{\zeta}^{\beta}  F_{Z_{E^c}} \Vert
 _{L^{\infty}(\R^{E}\times \R^E)} \ , \leqno (8.6)$$
where $\widetilde I_2(E)$ is the set of multi-indices $(\alpha ,
\beta)$ in $I_2(E)$, such that $\alpha_j + \beta_j \geq 1$ for all
$j\in E$. We also have, for any $0 \leq h < 1$ and for any  $F$ verifying the assumption $H_4(M, \varepsilon)$, the following estimates:
$$ \sup _{(\alpha , \beta)   \in
I_2(E) } h^{(|\alpha|+|\beta|)/2}
  \Vert
 \partial _z^{\alpha} \partial _{\zeta}^{\beta} T_h(E) F_{Z_{E^c}} \Vert
 _{L^{\infty}(\R^{E}\times \R^E)} \leq \sup _{(\alpha , \beta)\in \widetilde I_4(E)
 } h^{(|\alpha|+|\beta|)/2}
  \Vert
 \partial _z^{\alpha} \partial _{\zeta}^{\beta}  F_{Z_{E^c}} \Vert
 _{L^{\infty}(\R^{E}\times \R^E)}\ , 
 \leqno (8.7)$$
where $\widetilde I_4(E)$ is the set of every multi-indices $(\alpha
, \beta)$ in $I_4(E)$ satisfying $\alpha_j + \beta_j \geq 2$ for
each $j\in E$.

 }

\bigskip

{\it Proof.} We first note that:
$$ e^{{h\over 4}\Delta_j} -I ={h\over 4}\Delta_j V_{jh}$$
with the operator $V_{jh}$ given by:
$$(V_{jh} F)(x , \xi) = (\pi h)^{-1} \int _{\R^2\times [0, 1] } e^{-{1\over h} (
u^2 + v^2)} 2 \theta F ( x + \theta u e_j , \xi + \theta v e_j) du
dv d\theta\ , $$
where $(e_j)_{(j\in E)}$ is the canonical basis of $\R^E$.  Next, we observe that the operator $V_{jh}$ is bounded in
$L^{\infty}(\R^E \times \R^E)$ with a norm smaller than $1$.  Besides, the operator
$ e^{{h\over 4}\Delta_j} -I$ is also bounded in
$L^{\infty}(\R^E \times \R^E)$ with a norm smaller than $2$.  Consequently:
$$ \Vert
 \partial _z^{\alpha} \partial _{\zeta}^{\beta} T_h(E) F \Vert
 _{L^{\infty}(\R^{E}\times \R^E)} \leq 2^{|E|} \ \Vert (\prod _{j\in E} M_j) F
 \Vert _{L^{\infty}(\R^{E}\times \R^E)}  $$
for any multi-index $( \alpha , \beta)$  in
$I_2(E)$, where
 $$ M_j = {\partial^{\alpha_j} \over \partial z_j^{\alpha_j} }
 {\partial^{\beta_j} \over \partial \zeta _j^{\beta_j} } \ \ \ \ \
 {\rm if} \ \ \ \  \alpha _j + \beta _j \geq 1 \hskip 2cm M_j ={h\over 4}  \Delta _j
 \ \ \ \ \
 {\rm if}\ \ \ \  \alpha _j = \beta _j =0\ .$$
We then  deduce (8.6). When $F$ satisfies $H_4(M, \varepsilon)$, we also have:
$$ \Vert
 \partial _z^{\alpha} \partial _{\zeta}^{\beta} T_h(E) F \Vert
 _{L^{\infty}(\R^{E}\times \R^E)} \leq  \ \Vert (\prod _{j\in E} M'_j) F
 \Vert _{L^{\infty}(\R^{E}\times \R^E)}  $$
for every multi-index $( \alpha , \beta)$  in $I_2(E)$, with
 $$ M'_j ={h\over 4}  \Delta _j  {\partial^{\alpha_j} \over \partial z_j^{\alpha_j} }
 {\partial^{\beta_j} \over \partial \zeta _j^{\beta_j} }$$
 for all $j\in E$. Then we obtain (8.7) .

\bigskip

{\it Proof of Proposition 8.1.} We apply Proposition
8.3
with $G = T_h (E) F$ and Proposition 8.4 with $F$ where
  $F$ satisfies hypothesis $H_2(M , \varepsilon)$. We notice that the number of multi-indices $(\alpha , \beta)$ in $I_2(E)$ is
$9^{|E|}$. We also remark that, $h^{(|\alpha|+|\beta|)/2}
\leq h^{|E|/2}$, for $0 < h \leq 1$, for all  $(\alpha
, \beta)$  in $\widetilde I_2(E)$ and, $h^{(|\alpha|+|\beta|)/2} \leq h^{|E|}$,
for all  $(\alpha , \beta)$  in $\widetilde
I_4(E)$. We then observe that, under the assumption $H_2(M , \varepsilon)$, we have:
$$\Vert
 \partial _z^{\alpha} \partial _{\zeta}^{\beta}  F_{Z_{E^c}} \Vert
 _{L^{\infty}(\R^{E}\times \R^E)} \leq M K_2^{|E|} \prod _{j\in E}
 \varepsilon_j   \hskip 2cm (\alpha , \beta) \in \widetilde I_2(E)$$
and when $H_4(M , \varepsilon)$ is verified:
$$\Vert
 \partial _z^{\alpha} \partial _{\zeta}^{\beta}  F_{Z_{E^c}} \Vert
 _{L^{\infty}(\R^{E}\times \R^E)} \leq M K_4^{|E|} \prod _{j\in E}
 \varepsilon_j ^2  \hskip 2cm (\alpha , \beta) \in \widetilde I_4(E)\ .$$
Proposition 8.1 is then deduced.

\bigskip

{\it End of the proof of Proposition 7.2.} In view of Definition 5.1
concerning the hybrid quantization and proposition 8.1, together
with the Schur lemma, we deduce that:
$$ | < Op_h^{hyb, E} (T_h(E) F)f , g> |  \leq  M \left ( {C K_2 \sqrt {h}
\pi\over 2 } \right  )^{|E|} \left [ \prod _{j\in E}
\varepsilon_j\right ]   N_E(f) N_E(g) \leqno (8.8)$$
for every $f$ and $g$ in ${\cal
H}(\Gamma)$ and for all $F$ satisfying $H_2(M , \varepsilon)$,
where $C = 450$ is the constant appearing in Proposition 8.1, $K_2=
sup_{j\in \Gamma} \max (1 , \varepsilon_j^3)$,  and $N_E(f)$ is defined by (7.8). In the case when $F$ verifies $H_4(M , \varepsilon)$, we have:
$$ | < Op_h^{hyb, E} (T_h(E) F)f , g> |  \leq  M \left ( {C K_4 h
\pi\over 2 } \right  )^{|E|} \left [ \prod _{j\in E} \varepsilon_j
^2\right ]   N_E(f) N_E(g)\ , \leqno (8.9)$$
where $K_4 = sup_{j\in \Gamma}\max (1 , \varepsilon_j^6)$. When applying the Schur lemma, we have used the fact that:
$$(2 \pi h)^{-1/2} \int _{\R} \left ( 1 + {x^2 \over h} \right
)^{-1} dx = \sqrt {\pi \over 2}\ .$$
From Proposition 7.3, we have $N_E(f) = \Vert f \Vert _{{\cal
H}(\Gamma)} $ and Proposition 7.2 is proved.

\bigskip

{\bf 9. Comparison with previous definitions of the Weyl calculus.}

\bigskip

Some standard works consider the case when the symbol $F$ is a continuous function on $B(\Gamma)  \times B(\Gamma)$ and the Fourier transform of a bounded measure. We then assume that there exists a bounded measure   $\rho$ on $Z_{\R} \times Z_{\R}$,
where $Z_{\R} = \ell^2(\Gamma , \R)$, such that:
$$ F(x , \xi) = \int _{ Z_{\R} \times Z_{\R}}  e^{-i ( \ell_y ( x) + \ell_{\eta }(\xi))}
d\rho ( y , \eta) \ , \leqno (9.1)$$
where $\ell _y$ is the function defined on $B(\Gamma)$ in Theorem
2.6, for all $y$ in $Z_{\R}$. Thus, $F$ is a bounded function on
$B(\Gamma ) \times B(\Gamma)$.

\bigskip

When $E$ is a finite subset of  $\Gamma$, combining Theorem 3.4
(point (3.19)) and the definition (4.9) of the isomorphism
$\widetilde J_{Eh}^K$ concerning the Lebesgue measure, we have:
$$\left ( \widetilde J_{Eh}^K e^{{i\over \sqrt {h}} \Phi_S (a+i b)} \left (
 \widetilde J_{Eh}^K f \right ) ^{-1} \right ) (u) = f(u + b) e^{
 {i\over h} ( a \cdot u + {1 \over 2} a \cdot b) }$$
for any $a$ and $b$ in $\R^E$ and for every $f$
in $L^2(\R^E , \lambda_E)$.
Besides, it is well-known that the Weyl calculus has the following property (in finite dimension):
$$ E (x , \xi)= e^{ {i\over h}  (a \cdot x + b \cdot \xi)} \Longrightarrow
\left ( Op_h^{weyl}(E) f \right ) (u) =f(u + b) e^{
 {i\over h} ( a \cdot u + {1 \over 2} a \cdot b) }\ . $$
Therefore, it is natural to define the Weyl operator associated with
a symbol $F$ verifying (9.1) by:
$$ Op_h^{old-weyl}(F) = \int _{Z_{\R} \times Z_{\R}}
e^{-i\sqrt {h} \Phi_S (y + i \eta) } d\rho ( y , \eta)\ , \leqno (9.2)
$$
where   the   unbounded operator $\Phi _S ( y + i \eta)$, formally
self-adjoint, is associated with the element $y + i \eta$ of $Z_{\bf
C}$ as in (3.6), for all $(y , \eta)$ in $Z_{\R} \times Z_{\R}$. The
operator in  (9.2) is indeed bounded in the symmetric Fock space
 ${\cal H}(\Gamma)$ and we have:
$$ \Vert Op_h^{old-weyl}(F) \Vert _{{\cal L}({\cal H}(\Gamma))} \leq
 \int _{Z_{\R} \times
Z_{\R}} d|\rho | ( y , \eta)\ , \leqno (9.3)$$
where $|\rho|$ is the absolute value measure of the bounded measure
$\rho$. Definition (9.2) is considered by  Kree R\c akzka
[K-R], Lascar [LA1] and more recently by Albeverio Daletskii [A-D].

\bigskip

{\bf Theorem 9.1.} {\it Assume that  $F$ is a continuous function on
$B(\Gamma) \times B(\Gamma)$ written as in (9.1) (where $\rho$ is a
bounded measure on $Z_{\R} \times Z_{\R}$ with $Z_{\R} =
\ell^2(\Gamma , \R)$) and also verifying hypothesis $H_2(M,
\varepsilon)$ in Definition 1.3, where $(\varepsilon_j)_{(j\in
\Gamma)}$ is a summable family. Then, the two operators
$OP_h^{weyl}(F)$ and $OP_h^{old-weyl}(F)$, respectively defined by
(5.4) and (9.2), are equal.

}

\bigskip

Let  $(\Lambda_n)$ be an increasing sequence of finite subsets of $\Gamma$ with union $\Gamma$. For every $y$ in $\ell^2 (\Gamma , \R)$ and any $n\geq 0$, let:
$$ \left ( p_{n} (y) \right )_j =\left \{ \matrix { y_j &{\rm if}&
j \in \Lambda _n \cr  \cr  0 & {\rm if}& j\in \Lambda _n^c \cr
}\right .$$
and set $q_n = I - p_n$.

\bigskip

{\bf Lemma 9.2.} {\it Under the assumptions of Theorem 9.1, the
operator $Op_h^{hyb, \Lambda_n}(F)$ given by Definition 5.1 is
satisfying:
$$Op_h^{hyb , \Lambda_n}(F)  =  \int _{Z_{\R} \times Z_{\R}}
e^{-i\sqrt {h} \Phi_S (y + i \eta) } e^{ -{h\over 4} ( |q_n(y)|^2 +
|q_n(\eta )|^2 )} d\rho ( y , \eta)\ .\leqno (9.4)$$
}

\bigskip

{\it Proof of the lemma.}  We need to prove that:
$$ <Op_h^{hyb , \Lambda_n}(F)   f , g> =
\int _{Z_{\R} \times Z_{\R}}< e^{-i\sqrt {h} \Phi_S (y + i \eta) } f
, g> \ e^{ -{h\over 4} ( |q_n(y)|^2 + |q_n(\eta )|^2 )} d\rho ( y ,
\eta)\leqno (9.5)$$
for all $f$ and all
$g$ in ${\cal H}(\Gamma)$.
It suffices to prove this equality when:
$$ f = f _{\Lambda_n} \otimes f _{\Lambda_n^c } \hskip 2cm
g = g _{\Lambda_n} \otimes g _{\Lambda_n^c }\leqno (9.6)$$
with $f _{\Lambda_n}$ and $g _{\Lambda_n}$ in ${\cal H}(\Lambda
_n)$, $f _{\Lambda_n^c}$ and $g _{\Lambda_n^c}$ in ${\cal H}(\Lambda
_n^c)$. In this situation, we use Definition 5.1 of the hybrid
operator $Op_h^{hyb , \Lambda_n}(F)$. In this definition appears
the operator $Op_h^{weyl , \Lambda_n}(F_{Z_{\Lambda_n^c}})$. For
this operator, we may replace the definition in (5.1) by the one in
(9.2). These two definitions are indeed equivalent since $\Lambda_n$
is finite. We obtain:
$$<Op_h^{(hyb , \Lambda_n)}(F)   f , g> = \int _{B(\Lambda_n^c) \times B(\Lambda_n^c)
\times Z_{\R}\times Z_{\R}} <e^{-i\sqrt {h} \Phi_S (p_n (y + i \eta)
) }f _{\Lambda_n} , g _{\Lambda_n} > ...$$
$$ ...  \ e^{-i (\ell_{q_n( y)} (x_{\Lambda_n^c}) + \ell _{q_n ( \eta)}
(\xi  _{\Lambda_n^c} )) } \ \Big ( J^{\Phi} _{\Lambda _n , h} W
_{\Lambda_n^c}f _{\Lambda_n^c }\Big ) ( X_{\Lambda_n^c}) \ \overline
{ \Big ( J^{\Phi} _{\Lambda _n , h} W _{\Lambda_n^c}g _{\Lambda_n^c
}\Big ) ( X_{\Lambda_n^c})  } d\mu ^{\Phi}_{\Lambda_n^c , h} (
X_{\Lambda_n^c}) d\rho (y , \eta)\ .$$
From Theorem 4.8 applied with $E = \Lambda_n^c$ and $Y$ replaced by
$h q_n(Y)$, we see that:
$$e^{ -{h\over 4} ( |q_n(y)|^2 + |q_n(\eta )|^2 )}
<e^{-i\sqrt {h} \Phi_S (q_n (y + i \eta) ) }f _{\Lambda_n^c} , g
_{\Lambda_n^c} >= ... $$
$$ = \int _{B(\Lambda_n^c) \times B(\Lambda_n^c)}
 e^{-i (\ell_{q_n( y)} (x_{\Lambda_n^c}) + \ell _{q_n ( \eta)}
(\xi  _{\Lambda_n^c} )) }
 \ \Big ( J^{\Phi}
_{\Lambda _n , h} W _{\Lambda_n^c}f _{\Lambda_n^c }\Big ) (
X_{\Lambda_n^c}) \ \overline { \Big ( J^{\Phi} _{\Lambda _n , h} W
_{\Lambda_n^c}g _{\Lambda_n^c }\Big ) ( X_{\Lambda_n^c})  } d\mu
^{\Phi}_{\Lambda_n^c , h} ( X_{\Lambda_n^c})\ .
$$

Thus, we obtain (9.5) in the case (9.6). We then deduces (9.5) in the general case applying linearity, density and continuity arguments to  both sides.
\bigskip

{\it End of the proof of Theorem 9.1.} Suppose that
$Op_h^{old-weyl}(F)$ denotes the operator defined by the standard
relation (9.2) then, from (9.2), (9.4) and (9.3),
$$ \Vert Op_h^{(hyb , \Lambda_n)}(F) - Op_h ^{old-weyl}(F) \Vert _{{\cal L}({\cal H})} \leq
\int _{Z_{\R} \times Z_{\R} } \left | 1 - e^{ -{h\over 4} (
|q_n(y)|^2 + |q_n(\eta )|^2 )} \right | \  d|\rho | ( y , \eta)\ .$$
Lebesgue Theorem implies that
$$ \lim _{n\rightarrow +\infty}  \Vert Op_h^{(hyb , \Lambda_n)}(F) -
Op_h ^{old-weyl}(F) \Vert _{{\cal L}({\cal H}(\Gamma))} = 0\ . $$
In view of Theorem 5.4, we also have, since hypothesis $H_2(M,
\varepsilon)$ is satisfied and since $(\varepsilon_j)_{(j\in
\Gamma)}$ is a summable family,
$$ \lim _{n\rightarrow +\infty}  \Vert Op_h^{(hyb , \Lambda_n)}(F) -
Op_h ^{weyl}(F) \Vert _{{\cal L}({\cal H}(\Gamma))} = 0\ . $$
Therefore,  if  $F$ is as in (9.1) and verifies
hypothesis $H(M, \varepsilon)$, we deduce that the  operator $Op_h^{old-weyl}(F)$
given by the standard definition (which uses (9.1)) and  $Op_h^{weyl}(F)$ constructed in this work (which uses hypothesis $H(M, \varepsilon)$), are equal.

\bigskip

{\bf References.}

\bigskip

[A-D] S. Albeverio, A. Daletskii, {\it Algebras of
pseudodifferential operators in $L^2$ given by smooth measures on
Hilbert spaces,} Math. Nachr. {\bf 192} (1998) 5-22.

\smallskip

[A] Z. Ammari, {\it Canonical commutation relations and interacting
Fock spaces,} Journ\'ees E.D.P, exp. 2, 13p, Ecole Polytechnique,
Palaiseau, 2004.

\smallskip

[A-B] Z. Ammari, S. Breteaux, {\it Propagation of chaos for
many-bosons systems in one dimension with a point-pair interaction,}
arXiv:0906.3047v1 (2009).

\smallskip

[B-Z-H] J. Baez, I. E. Segal, Z-F Zhou, {\it Introduction to
algebraic and constructive quantum field theory},  Princeton Series
in Physics. Princeton University Press, Princeton, NJ, 1992.

\smallskip

[B-K] Y. M. Berezansky, Y. G. Kondratiev, {\it  Spectral methods in
infinite-dimensional analysis. Vol. 1 and Vol 2,} Translated from
the 1988 Russian original by P. V. Malyshev and D. V. Malyshev.
Mathematical Physics and Applied Mathematics, 12/2. Kluwer Academic
Publishers, Dordrecht, 1995.

\smallskip

[C-V] A.P. Calder\'on, R. Vaillancourt, {\it A class of bounded
pseudo-differential operators,} Proc. Nat. Acad. Sci. U.S.A, {\bf
69}, (1972), 1185-1187.

\smallskip

[C-M] R. L. Coifman, Y. Meyer, {\it Au del\`a des op\'erateurs
pseudo-diff\'erentiels,} Ast\'erisque, {\bf 57}, 1978.

\smallskip

[C] H. O. Cordes, {\it On compactness of commutators of
multiplications and convolutions, and boundedness of
pseudo-differential operators,} J. Funct. Anal, {\bf 18} (1975)
115-131.

\smallskip

[D-Z] J. Derezinski, C. G\'erard,{\it  Asymptotic completeness in
quantum field theory. Massive Pauli-Fierz Hamiltonians.}  Rev. Math.
Phys. {\bf 11} (1999), no. 4, 383-450.

\smallskip

[D-H] B. Driver, B. Hall, {\it Yang-Mills theory and the Segal
Bargmann  transform,} Comm. in Math. Phys, {\bf 201} (1999) (2)
249-290.

\smallskip

[FA] J. Faraut, {\it Espaces Hilbertiens invariants de fonctions
holomorphes,} S\'eminaires \& Congr\`es, {\bf 7}, (2003), 101-167.

\smallskip

[FO]G. B. Folland, {\it Harmonic Analysis on phase space,} Annals of
Mathematics studies {\bf 122}, Princeton University Press, Princeton
(N.J), 1989.

\smallskip

[G1] L. Gross {\it Measurable functions on Hilbert space, } Trans.
Amer. Math. Soc. {\bf 105} (1962) 372–390.

\smallskip

[G2] L. Gross, {\it Abstract Wiener spaces}, Proc. 5th Berkeley Sym.
Math. Stat. Prob, {\bf 2}, (1965), 31-42.

\smallskip

[G3]L. Gross,  {\it Abstract Wiener measure and infinite dimensional
potential theory}, in {\it Lectures in modern Analysis and
applications, II,} Lecture Notes in Math {\bf 140}, 84-116, Springer
(1970).

\smallskip

[HA] B. Hall, {\it  Holomorphic methods in analysis and mathematical
physics,} First Summer School in Analysis and Mathematical Physics
(Cuernavaca Morelos, 1998), 1–59, Contemp. Math., {\bf 260}, Amer.
Math. Soc., Providence, RI, 2000.

\smallskip

[HO] L. H\"ormander, {\it The analysis of linear partial
differential operators,} Volule III, Springer, 1985.

\smallskip

[HW] I. L. Hwang, {\it The $L^2$ boundedness of pseudo-differential
operators,} Trans. Amer. Math. Soc, {\bf 302} (1987) 55-76.

\smallskip

[J] S. Janson, {\it Gaussian Hilbert spaces,} Cambridge Tracts in
Math. {\bf 129}, Cambridge Univ. Press (1997).

\smallskip

[KH] A. Khrennikov, {\it Distributions and pseudo-differential
operators on infinite-dimensional spaces with applications in
quantum physics,} in  Pseudo-differential operators and related
topics, 161-172, Oper. Theory Adv. Appl., {\bf 164}, Birkh\"auser,
Basel, 2006.

\smallskip

[K-R] P. Kr\'ee, R.  R\c aczka, {\it  Kernels and symbols of
operators in quantum field theory,} Ann. Inst. H. Poincar\'e Sect. A
(N.S.) {\bf 28} (1978), no. 1, 41–73.

\smallskip

[KU] H. H. Kuo, {\it Gaussian measures in Banach spaces.}  Lecture
Notes in Mathematics, Vol. 463. Springer, Berlin-New York, 1975.

\smallskip

[LA1] B. Lascar, {\it Noyaux d'une classe d'op\'erateurs
pseudo-diff\'erentiels sur l'espace de Fock, et applications.}
S\'eminaire Paul Kr\'ee, 3e ann\'ee (1976-77), Equations aux
d\'eriv\'ees partielles en dimension infinie, Exp. No. 6, 43 pp.

\smallskip

[LA2] B. Lascar,  {\it Equations aux d\'eriv\'ees partielles en
dimension infinie.}  Vector space measures and applications (Proc.
Conf., Univ. Dublin, Dublin, 1977), I, pp. 286-313, Lecture Notes in
Math, {\bf 644}, Springer, Berlin, 1978.

\smallskip

[LA3] B. Lascar,  {\it Op\'erateurs pseudo-diff\'erentiels en
dimension infinie. Etude de l'hypoellipticit\'e de la
r\'esolubilit\'e dans des classes de fonctions holderiennes et de
distributions pour des op\'erateurs pseudo-diff\'e\-rentiels
elliptiques,}  J. Analyse Math. {\bf 33} (1978), 39-104.

\smallskip

[LA4] B. Lascar,  {\it Une classe d'op\'erateurs elliptiques du
second ordre sur un espace de Hilbert,}  J. Funct. Anal. {\bf 35}
(1980), no. 3, 316-343.

\smallskip

[LA5] B. Lascar, {\it  Probl\`emes de Cauchy hyperboliques en
dimension infinie,}  S\'eminaire Paul Kr\'ee, 4e ann\'ee: 1977-1978.
Equations aux d\'eriv\'ees partielles en dimension infinie, Exp. No.
6, 29 pp.

\smallskip

[LA6] B. Lascar, {\it Th\'eor\`eme de Cauchy-Kovalevsky et
th\'eor\`eme d'unicit\'e d'Holmgren pour des fonctions analytiques
d'une infinit\'e de variables.}   in holomorphy (Proc. Sem. Univ.
Fed. Rio de Janeiro, Rio de Janeiro, 1977), pp. 485-508,
North-Holland Math. Stud., {\bf 34}, North-Holland, Amsterdam-New
York, 1979.

\smallskip

[LA7] B.  Lascar, {\it  Invariance par diff\'eomorphisme d'espace de
Sobolev. Espace de Sobolev d'une vari\'et\'e. Applications. }
S\'eminaire Paul Kr\'ee 2e ann\'ee (1975-76), Equations aux
d\'eriv\'ees partielles en dimension infinie, Exp. No. 7, 29 pp.

\smallskip

[LA8] B.  Lascar, {\it Une condition n\'ecessaire et suffisante
d'ellipticit\'e pour une classe d'op\'erateurs diff\'erentiels en
dimension infinie.}  Comm. Partial Differential Equations {\bf 2}
(1977), no. 1, 31-67.

\smallskip

[LA9] B. Lascar, {\it Th\'eor\`eme de Cauchy-Kovalevsky et
th\'eor\`eme d'unicit\'e d'Holmgren pour des fonctions analytiques
d'une infinit\'e de variables.}  S\'eminaire Paul Kr\'ee 2e ann\'ee
(1975-76), Equations aux d\'eriv\'ees partielles en dimension
infinie, Exp. No. 8, 16 pp.

\smallskip

[LA10] B.  Lascar, {\it  Op\'erateurs pseudo-diff\'erentiels en
dimension infinie. Applications. } C. R. Acad. Sci. Paris S\'er. A-B
{\bf 284} (1977), no. 13, A767-A769,

\smallskip

[LA11] B. Lascar, {\it  M\'ethodes $L^{2}$ pour des \'equations aux
d\'eriv\'ees partielles d\'ependant d'une infinit\'e de variables.}
 S\'eminaire Goulaouic-Schwartz (1975-1976) Equations aux
d\'eriv\'ees partielles et analyse fonctionnelle, Exp. No. 5, 11 pp.
Centre Math., Ecole Polytech., Palaiseau, 1976.

\smallskip

[LA12] B.  Lascar, {\it  Propri\'et\'es locales d'espaces de type
Sobolev en dimension infinie. } Comm. Partial Differential Equations
{\bf 1} (1976), no. 6, 561-584.

\smallskip

[LA13] B. Lascar, {\it  Th\'eor\`eme de Cauchy-Kovalewsky et
th\'eor\`eme d'unicit\'e d'Holmgren pour des fonctions d'une
infinit\'e de variables.}  C. R. Acad. Sci. Paris S\'er. A-B {\bf
282} (1976), no. 13, A691-A694.

\smallskip

[LA14] B. Lascar, {\it  Propri\'et\'es locales d'espaces de type
Sobolev en dimension infinie. } S\'eminaire Paul Kr\'ee, 1re ann\'ee
(1974-75), Equations aux d\'eriv\'ees partielles en dimension
infinie, Exp. No. 11, 16 pp.

\smallskip

[LA15] B. Lascar, {\it  Op\'erateurs pseudo-diff\'erentiels d'une
infinit\'e de variables, d'apr\`es M. I. Visik.} S\'eminaire Pierre
Lelong (Analyse), Ann\'ee 1973-1974, pp. 83–90. Lecture Notes in
Math., {\bf 474}, Springer, Berlin, 1975.

\smallskip

[LA16] B. Lascar, {\it  Propri\'et\'es d'espaces de Sobolev en
dimension infinie, } C. R. Acad. Sci. Paris S\'er. A-B {\bf 280}
(1975), no. 23, A1587-A1590.

\smallskip

[LER] N. Lerner, {\it Metrics on the phase space and non
self-adjoint pseudo-differential operators,} Birkh\"auser Springer,
2010.

\smallskip

[LEV] Th. L\'evy, {\it Mesures gaussiennes et espaces de Fock}
(Peyresq, 2003, preprint).

\smallskip

[RE-SI] M. Reed, B. Simon, {\it Methods of modern mathematical
physics, Vol II, Fourier Analysis, self-adjointness,} Academic
Press, 1975.

\smallskip

[RO-SC] I. Rodnianski, B. Schlein, {\it Quantum fluctuations and
rate of convergence towards mean field dynamics,} Comm. Math. Phys,
{\bf 291} (2009) 31-61.

\smallskip

[SE] I. Segal, {\it Tensor algebras over Hilbert spaces, I,} Trans.
Amer. Math. Soc, {\bf 81} (1956), 104-134.

\smallskip

[SI1] B. Simon, {\it Functional integration and quantum physics,}
Second Edition, AMS Chelsea Publ, Providence (R.I), 2005.

\smallskip

[SI2] B. Simon, {\it The $P(\varphi)_2$ Euclidean quantum field
theory}, Princeton series in Physics, Princeton Univ. Press,
Princeton (N.J) 1974.

\smallskip

[SJ] J. Sj\"ostrand, {\it Singularit\'es analytiques microlocales,}
Ast\'erisque {\bf 95} (1982).

\smallskip

[U1] A. Unterberger, {\it  Oscillateur harmonique et op\'erateurs
pseudo-diff\'erentiels},   Ann. Inst. Fourier (Grenoble) {\bf 29}
(1979), no. 3, xi, 201–221.

\smallskip

[U2] A. Unterberger, {\it Les op\'erateurs m\'etadiff\'erentiels},
in Complex analysis, microlocal calculus and relativistic quantum
theory, Lecture Notes in Physics {\bf 126} (1980) 205-241.

\bigskip

Laboratoire de Math\'ematiques, FR CNRS 3399, EA 4535, Universit\'e de Reims
Champagne-Ardenne, Moulin de la Housse, B. P. 1039, F-51687
Reims, France, 

{\it E-mail:} {\tt laurent.amour@univ-reims.fr }

{\it E-mail:} {\tt lisette.jager@univ-reims.fr}

{\it E-mail:} {\tt jean.nourrigat@univ-reims.fr}

\end